\documentclass[10pt]{article}
\date{}

\usepackage{hyperref}
\usepackage{amssymb, amsmath}
\newcommand{\Pf}{\mathop{{\rm Pf}}\nolimits}

\newcommand{\cDt}{\widetilde{\cD}}
\newcommand{\mlabel}{\label}

\usepackage{amsmath}
\usepackage{latexsym}
\usepackage{amssymb}

\usepackage{amsthm}

\font\tengoth=eufm10 at 10pt
\font\sevengoth=eufm7 at 6pt
\newfam\gothfam
\textfont\gothfam=\tengoth
\scriptfont\gothfam=\sevengoth

\newcommand{\g}{{\mathfrak g}}

\newcommand{\h}{{\mathfrak h}}

\newcommand{\fg}{{\mathfrak g}}
\newcommand{\fh}{{\mathfrak h}}

\newcommand{\fo}{{\mathfrak o}}
\newcommand{\fq}{{\mathfrak q}}

\renewcommand{\:}{\colon}
\newcommand{\1}{\mathbf{1}}

\newcommand{\cA}{\mathcal{A}}
\newcommand{\cB}{\mathcal{B}}

\newcommand{\cD}{\mathcal{D}}
\newcommand{\cE}{\mathcal{E}}
\newcommand{\cF}{\mathcal{F}}

\newcommand{\cH}{\mathcal{H}}

\newcommand{\cL}{\mathcal{L}}
\newcommand{\cM}{\mathcal{M}}
\newcommand{\cN}{\mathcal{N}}
\newcommand{\cO}{\mathcal{O}}
\newcommand{\cP}{\mathcal{P}}

\newcommand{\cS}{\mathcal{S}}

\newcommand{\cW}{\mathcal{W}}

\newcommand{\eset}{\emptyset}

\renewcommand{\phi}{\varphi}
\newcommand{\dd}{{\tt d}}

\newcommand{\subeq}{\subseteq}
\newcommand{\supeq}{\supseteq}

\newcommand{\into}{\hookrightarrow}
\newcommand{\eps}{\varepsilon}

\newcommand{\N}{{\mathbb N}}
\newcommand{\Z}{{\mathbb Z}}
\newcommand{\R}{{\mathbb R}}
\newcommand{\C}{{\mathbb C}}

\newcommand{\bP}{{\mathbb P}}

\newcommand{\bS}{{\mathbb S}}

\renewcommand{\hat}{\widehat}

\renewcommand{\tilde}{\widetilde}

\renewcommand{\L}{\mathop{\bf L{}}\nolimits}

\newcommand{\Aff}{\mathop{{\rm Aff}}\nolimits}

\newcommand{\GL}{\mathop{{\rm GL}}\nolimits}
\newcommand{\SL}{\mathop{{\rm SL}}\nolimits}

\newcommand{\SO}{\mathop{{\rm SO}}\nolimits}
\newcommand{\SU}{\mathop{{\rm SU}}\nolimits}
\newcommand{\OO}{\mathop{\rm O{}}\nolimits}

\newcommand{\U}{\mathop{\rm U{}}\nolimits}

\newcommand{\Fix}{\mathop{{\rm Fix}}\nolimits}

\newcommand{\Ad}{\mathop{{\rm Ad}}\nolimits}

\newcommand{\Herm}{\mathop{{\rm Herm}}\nolimits}

\newcommand{\Diff}{\mathop{{\rm Diff}}\nolimits}

\newcommand{\diag}{\mathop{{\rm diag}}\nolimits}

\newcommand{\id}{\mathop{{\rm id}}\nolimits}

\newcommand{\supp}{\mathop{{\rm supp}}\nolimits}

\newcommand{\Spann}{\mathop{{\rm span}}\nolimits}
\newcommand{\ev}{\mathop{{\rm ev}}\nolimits}

\newcommand{\oline}{\overline}

\newcommand{\la}{\langle}
\newcommand{\ra}{\rangle}

\newcommand{\res}{\vert}

\newcommand{\ssssarr}{\hbox to 15pt{\rightarrowfill}}
\newcommand{\sssarr}{\hbox to 20pt{\rightarrowfill}}
\newcommand{\ssarr}{\hbox to 30pt{\rightarrowfill}}
\newcommand{\sarr}{\hbox to 40pt{\rightarrowfill}}
\newcommand{\arr}{\hbox to 60pt{\rightarrowfill}}
\newcommand{\larr}{\hbox to 60pt{\leftarrowfill}}
\newcommand{\Arr}{\hbox to 80pt{\rightarrowfill}}

\def\theoremname{Theorem}
\def\propositionname{Proposition}
\def\corollaryname{Corollary}
\def\lemmaname{Lemma}
\def\remarkname{Remark}
\def\conjecturename{Conjecture} 

\def\definitionname{Definition}
\def\exercisename{Exercise}
\def\examplename{Example}
\def\examplesname{Examples}
\def\problemname{Problem}
\def\problemsname{Problems}
\def\proofname{Proof}

\def\satzname{Satz} 
\def\koroname{Korollar}
\def\folgname{Folgerung}
\def\bemerkname{Bemerkung}
\def\aufgname{Aufgabe}

\def\beisname{Beispiel}
\def\beissname{Beispiele}
\def\bewname{Beweis}

\def\@thmcounter#1{\noexpand\arabic{#1}}
\def\@thmcountersep{}
\def\@begintheorem#1#2{\it \trivlist \item[\hskip 
\labelsep{\bf #1\ #2.\quad}]}
\def\@opargbegintheorem#1#2#3{\it \trivlist
      \item[\hskip \labelsep{\bf #1\ #2.\quad{\rm #3}}]}
\makeatother
\newtheorem{theor}{\theoremname}[section]
\newtheorem{propo}[theor]{\propositionname}
\newtheorem{coro}[theor]{\corollaryname}
\newtheorem{lemm}[theor]{\lemmaname}

\newenvironment{thm}{\begin{theor}\it}{\end{theor}}

\newenvironment{prop}{\begin{propo}\it}{\end{propo}}

\newenvironment{cor}{\begin{coro}\it}{\end{coro}}

\newenvironment{lem}{\begin{lemm}\it}{\end{lemm}}

\newtheorem{rema}[theor]{\remarkname}

\newenvironment{rem}{\begin{rema}\rm}{\end{rema}}

\newtheorem{stepnow}[theor]{}

\newtheorem{defin}[theor]{\definitionname} 

\newenvironment{defn}{\begin{defin}\rm}{\end{defin}}

\newtheorem{exerc}{\exercisename}[section]

\newtheorem{exa}[theor]{\examplename}

\newenvironment{ex}{\begin{exa}\rm}{\end{exa}}

\newtheorem{exas}[theor]{\examplesname}

\newenvironment{exs}{\begin{exas}\rm}{\end{exas}}

\newtheorem{conj}[theor]{\conjecturename}

\newtheorem{pro}[theor]{\problemname}

\newtheorem{prs}[theor]{\problemsname}

\newenvironment{Proof*}{\begin{trivlist}\item[\hskip%
\labelsep{\bf\proofname.\quad}]}%
{\end{trivlist}}

\newenvironment{prf}{\begin{proof}}{\end{proof}}
  
\newcommand{\pmat}[1]{\begin{pmatrix} #1 \end{pmatrix}}

{\hfill\qed\end{trivlist}}

\newenvironment{beweis*}{\begin{trivlist}\item[\hskip%
\labelsep{\bf\bewname.\quad}]}%
{\end{trivlist}}

\newtheorem{satzn}[theor]{\satzname}

\newtheorem{koro}[theor]{\koroname}

\newtheorem{folg}[theor]{\folgname}

\newtheorem{bem}[theor]{\bemerkname}

\newtheorem{aufg}[theor]{\aufgname}

\newtheorem{aufgn}[theor]{\aufgname}

\newtheorem{beis}[theor]{\beisname}

\newtheorem{beiss}[theor]{\beissname}

\begin{document}

\title{Reflection Positivity and Conformal Symmetry}
\author{Karl-Hermann Neeb,
\begin{footnote}{
Department  Mathematik, FAU Erlangen-N\"urnberg, Cauerstrasse 11, 
91058-Erlangen, Germany; neeb@math.fau.de}
\end{footnote}
\begin{footnote}{Supported by DFG-grant NE 413/7-2, Schwerpunktprogramm 
``Darstellungstheorie''.} 
\end{footnote}
Gestur \'Olafsson
\begin{footnote}{Department of mathematics, Louisiana State University, 
Baton Rouge, LA 70803, USA; olafsson@math.lsu.edu}
\end{footnote}
\begin{footnote}
{The research of G. \'Olafsson was supported by NSF grants 
DMS-0801010, DMS-110337 and the Emerging Fields Project 
``Quantum Geometry'' of the University of Erlangen.} 
\end{footnote}
}

\maketitle

\abstract{A reflection positive Hilbert space is a triple $(\cE,\cE_+,\theta)$, where 
$\cE$ is a Hilbert space, $\theta$ a unitary involution and 
$\cE_+$ a closed subspace on which the hermitian form 
$\la v,w\ra_\theta := \la \theta v, w \ra$ is positive semidefinite. From this data one 
obtains a  Hilbert space $\hat\cE$ by completing a suitable quotient of 
$\cE_+$ with respect to $\la \cdot,\cdot \ra_\theta$ on 
$\cE_+$. To obtain compatible unitary representations of Lie groups, we start with triples 
$(G,S,\tau)$, where $G$ is a Lie group, $\tau$ an involutive automorphism of $G$ 
and $S$ a subsemigroup invariant under the involution $s^\sharp = \tau(s)^{-1}$. 
Then a unitary representation $\pi$ of $G$ on $(\cE,\cE_+,\theta)$ is called reflection 
positive if $\theta \pi(g) \theta= \pi(\tau(g))$ and $\pi(S)\cE_+ \subeq \cE_+$. 
Motivated by the passage from the euclidean motion group to the Poincar\'e group 
in quantum field theory, one expects a duality between reflection positive representations 
and unitary representations of the dual symmetric Lie group $G^c$ on $\hat\cE$. 

We propose a new approach to a reflection positive 
representations based on reflection positive distributions and 
reflection positive distribution vectors. 
In particular, we generalize the Bochner--Schwartz Theorem 
to positive definite distributions on open convex cones and apply our 
techniques to complementary series representations of the conformal group $\OO_{1,n+1}^+(\R)$ 
of the sphere $\bS^n$.}

\section*{Introduction} \mlabel{sec:0}

The concept of reflection positivity has its origins in the work of Osterwalder--Schrader \cite{OS73, OS75}  
on constructive quantum field theory and duality
between unitary representations of the euclidean motion group $\cE_{n}=\OO_n (\R)\ltimes \R^n$ 
and the Poincar\'e group $P_n=\OO^+_{1,n-1}(\R) \ltimes\R^{1,n-1}$ of affine isometries of 
$n$-dimensional Minkowski space. 
Here $\OO^+_{1,n-1}(\R) $ is the group preserving the Lorentz form $(t,x) \mapsto t^2-\|x\|^2$ and mapping the 
forward light cone \[ \Omega =\{(t,x)\mid t^2-\|x\|^2>0, t>0\}\] onto itself. 
Multiplying the time coordinate $t$ by $i$ transform
the Lorentz form into  $-t^2-\|x\|^2=-\|(t,x)\|^2$ setting up a duality between the groups 
$P_n$ and $E_n$.

On the mathematical side this duality can be made precise as follows. 
If $\g$ is a Lie algebra with an involutive automorphism $\tau$, then 
we have the $\tau$-eigenspace decomposition 
$\g = \fh \oplus \fq = \ker(\tau - \1) \oplus \ker(\tau + \1)$ and the subspace 
$\g^c := \fh \oplus i \fq$ of the complexification $\g_\C$ of $\g$ is another real form. 
We thus obtain a duality between the pairs $(\g,\tau)$ and $(\g^c,\tau)$. 
At the core of the notion of reflection positivity is the idea that 
this duality can sometimes be implemented on the level of unitary representations. 
This is quite simple on the Lie algebra level: 
Let $\cE^0$ be a pre-Hilbert space and $\pi$ be a representation 
of $\g$ on $\cE^0$ by skew-symmetric operator. We also assume that 
there exists a unitary operator $\theta$ on $\cE^0$ with 
$\theta\pi(x) \theta = \pi(\tau x)$ for $x \in \g$, and a $\g$-invariant subspace 
$\cE^0_+$ which is {\it $\theta$-positive} in the sense that the hermitian form  
$\la v, w \ra_\theta := \la \theta v,w\ra$ is positive semidefinite on $\cE^0_+$. Then complex 
linear extension leads to a representation of $\g^c$ on $\cE^0_+$ by 
operators which are skew-symmetric with respect to $\la \cdot, \cdot \ra_\theta$, so that we obtain a 
``unitary'' representation of $\g^c$ on the pre-Hilbert space $\hat\cE^0 := \cE^0_+/\{ v \: 
\cE^0_+ \: \la \theta v, v \ra =0\}$. 
This is the basic idea behind the reflection positivity correspondence between 
unitary representations of $\g$ and $\g^c$. What this simple picture completely ignores are 
issues of integrability and essential selfadjointness of operators. There are various 
natural ways to address these problems. Important first steps in this direction have been 
undertaken by Klein and Landau in \cite{KL81, KL82}, and 
Fr\"ohlich, Osterwalder and Seiler introduced in \cite{FOS83}
 the concept of a virtual representation, 
which was developed in greated generality  by Jorgensen in \cite{Jo86, Jo87}. 

The approach we shall pursue in the present paper is based on open subsemigroups 
of Lie groups. Starting with a unitary representation 
$(\pi,\cE)$ of a Lie group $G$ and a unitary involution $\theta$ of $\cE$ 
with $\theta\pi(g) \theta = \pi(\tau(g))$ for an involutive automorphism $\tau$ of $G$, 
we are looking for $\theta$-positive closed subspaces $\cE_+ \subeq \cE$ 
that are invariant under a subsemigroup $S \subeq G$ which is invariant 
under the involution $s \mapsto s^\sharp := \tau(s)^{-1}$. This leads to a 
$\sharp$-representation of $S$ by contractions on the Hilbert space $\hat\cE$, and there are powerful 
tools based on the L\"uscher--Mack Theorem to derive from such representations unitary representations 
of a Lie group $G^c$ with Lie algebra $\g^c$ 
(cf.\ \cite{LM75}, \cite[Sect.~9.5]{HN93} and \cite{MN11} for a generalization to Banach--Lie groups). 
We therefore focus on the triple $(G,\tau, S)$ and unitary representations as above.
In any case, one obtains unitary representations $\pi^c$ of $G^c$ for which 
the operators $-i\dd\pi^c(x)$ for $x \in i \fq$ and $\exp(\R_+ i x) \subeq S$ 
have positive spectrum. This condition imposes serious restrictions on 
the unitary representations of $G^c$ that one can obtain in this context. 

One of the basic requirements of quantum theory is that physical states form a 
Hilbert space $\cH$ with a unitary positive energy representation
$\pi$ of the Poincar\'e group, i.e., the spectral measure 
of the translation group is supported by the future light cone. 
Here the involution $\tau$ corresponds to time reversal 
which implements the duality between the Lie algebras of $E_n$ and $P_n$. 
Therefore one is lead to subsemigroups $S \subeq E_n$ containing the 
ray $\{ (x_0, 0) \in \R^n\: x_0 > 0\}$ which under $c$-duality corresponds to time translations. 
A natural enlargement is the semigroup 
\[ S=\OO_{n-1}(\R) \ltimes \{(t,x)\mid t>0,x\in\R^{n-1}\},\] 
but here the semigroup approach is problematic because there is no proper semigroup 
with interior points in $E_n$. This makes the passage from unitary representations  
of $E_n$ to unitary representation of $P_n$ a harder problem which can be addressed with the 
methods developed in \cite{Jo86, Jo87}. 
For further reference on the physical side of 
reflection positivity we would like to point out the work of A.~Klein 
\cite{Kl77,Kl78}, 
as well as the more resent  work by Jaffe and Ritter \cite{JR08,JR07}. An excellent
introduction can be found in the overview article \cite{JA08}, in particular Section~VII.
The approach to reflection positivity in terms of $c$-duality of Lie algebras and contractive 
representations of semigroups was already taken up in the work by Schrader in \cite{Sch86}, 
where he used reflection positivity to construct from a complementary series representation of 
$\SL_{2n}(\C)$ a unitary representation of the $c$-dual group 
$\SU_{n,n}(\C) \times \SU_{n,n}(\C)$ but without identifying 
the resulting representation. This approach was developed more systematically 
in \cite{JO98,JO00} where duality between semisimple causal symmetric spaces, the theory of compression
semigroups and the L\"uscher--Mack Theorem were used to construct from a generalized 
complementary series representation related to an ordered symmetric space $G/H$ 
an irreducible unitary highest weight (=positive energy) representation
of the dual group $G^c$. The construction was carried out for the case where 
$G$ is locally isomorphic to the automorphism group of a tube domain $T_\Omega := \R^k+i\Omega$, 
where $\Omega \subeq \R^k$ is a symmetric cone and $H= \{g\in G\mid g(\Omega)=\Omega\}$ 
is the subgroup preserving the totally real submanifold $\Omega$ of $T_\Omega$. In this case 
$\g^c \cong \g$, so that one obtains a correspondence between unitary representations 
of the same group. 

We have already mentioned that there are spectral restrictions on the representations 
of $G^c$ that can result from reflection positivity and this also 
restricts the class of groups for which this process can possibly apply. 
For semisimple Lie groups an inspection of parameters of those representations, 
in particular the infinitesimal character, indicates that the representations 
of $G$ to start with should be generalized complementary series representations. 
This partly explains why we as well as \cite{Sch86} and \cite{JO98,JO00} 
start with this class of representations. 

In the present article we first propose a new approach to a 
systematic treatment of reflection positive representations of triples 
$(G,S,\tau)$ based on {\it reflection positive distributions} 
and {\it reflection positive distribution vectors} of a unitary representation of $G$. 
Here the key tool is a refinement of the well-known GNS construction for positive definite 
functions and distributions to the reflection positive setting. To understand the 
nature of reflection positive distributions, it is indispensible to have a complete 
picture of the abelian case. To this end we obtain integral representations of 
reflection positive functions on the real line and generalize the Bochner--Schwartz Theorem 
to positive definite distributions on open convex cones. For the verification of positive 
definiteness of holomorphic kernels, we provide in Appendix~\ref{app:b} an general theorem 
asserting that it suffices to verify positive definiteness on open subsets or 
on totally real submanifolds. Finally we apply these techniques to exhibit some of 
the complementary series representations of the conformal group $\OO_{1,n+1}^+(\R)$ 
of the sphere $\bS^n$ as reflection positive. 

This article is organized as follows. 
In Section~\ref{sec:1} we  introduce reflection positive unitary representations 
for triples $(G,S,\tau)$ consisting of a Lie group $G$, 
an involutive automorphism $\tau$ and a subsemigroup $S$ invariant under the involution 
$s^\sharp := \tau(s)^{-1}$. To develop a systematic approach to reflection positivity 
for unitary representations, we introduce a concept of {\it reflection positive} 
operator-valued positive definite functions on $G$  and 
explain how they correspond to reflection positive representations $(\pi, \cE)$ 
that are generated in a very controlled fashion by a subspace $\cF$ of $\theta$-fixed 
vectors. This reflection positive 
variant of the GNS construction Proposition~\ref{prop:1.xx} is the main point 
of Section~\ref{sec:1}. 
In many interesting situations (see Section~\ref{sec:5}) 
it turns out that the generating vectors cannot be found in the Hilbert space itself 
but have to be replaced by distribution vectors in the larger space 
$\cE^{-\infty}$. This leads us in Section~\ref{sec:2} to the notion of a 
{\it reflection positive distribution}. On the representation side, 
they correspond to {\it reflection positive distribution cyclic representations}, 
which are triples $(\pi, \cE, \alpha)$, where $\pi$ is a unitary representation of 
$G$ satisfying $\theta \pi(g) \theta = \pi(\tau(g))$ for a unitary involution 
$\theta$ on $\cE$, 
$\alpha\in \cE^{-\infty}$ is a $\theta$-invariant cyclic distribution vector, 
and the closed subspace $\cE_+$ generated by $\pi^{-\infty}(\cD(S))\alpha$ is 
$\theta$-positive. Here $\cD(S) =C^\infty_c(S)$ denotes the space of test functions on $S$. 
In Proposition \ref{prop:1.y}, the main result of Section~\ref{sec:2}, 
we show that distribution 
cyclic representations are in one-to-one correspondence with reflection positive 
distributions on $G$. 

Section~\ref{sec:3} and \ref{app:a}  are devoted to classifying reflection positive functions 
in the finite-dimensional abelian case. For the triple $(\R,-\id,\R_+)$ we obtain 
complete information on reflection positive functions in terms of an integral representation, 
in which the building blocks are the functions $\phi(x) = e^{-\lambda|x|}$, $\lambda \geq 0$, 
for which the associated $\R_+$-representation on $\hat\cE$ is one-dimensional 
(Proposition~\ref{prop:2.1b}). This generalizes results of A.~Klein \cite{Kl77} obtained 
in the context of (OS)-positive covariance functions. 
We connect this case to \cite{JO00} and the reflection positivity on the group 
$\SL_2(\R)$ by considering natural abelian subgroups. One should also note that the 
$ax+b$-group, i.e., the affine group of the real line, is a subgroup of $\SL_2(\R)$, 
so that \cite{JO00} also gives rise to reflection positivity for this group.

The classification of reflection positive distributions is rather subtle. 
Here one can only hope to get hold of the restrictions to $S$ and, already for 
the open half line $S = \R_+  \subeq \R = G$, the classification of reflection 
positive extensions of positive definite distributions on $S$ to $G$ 
seems to be a hopeless task. On the other hand, the restriction to $S$ carries all 
information required for the representation of $G^c$, so that one is rather interested in 
``natural'' extensions of distributions from $S$ to $G$ and not in all of them. 
This motivates the main result of Section~\ref{app:a} 
which is a generalization of the Bochner--Schwartz Theorem (Theorem~\ref{thm:a.7}). 
For the case where $G = V$ is a vector space and 
$S = \Omega$ is an open convex cone invariant under $s^\sharp = -\tau(s)$, 
it provides an integral representation of positive definite distributions on $S$ 
that lead to contraction representations of $S$. This is precisely the class of distributions 
showing up for the representations on the spaces $\hat\cE$. For $\tau = -\id_V$ and open 
cones not containing affine lines we obtain prove the considerably stronger result that 
positive definite distributions are actually analytic functions (Theorem~\ref{thm:4.1}). 

In the last two sections we discuss reflection positive distribution vectors 
for the complementary series representations 
of the conformal group of $\R^n$, resp., its conformal compactification $\bS^n$. 
Here we start from canonical kernel functions on $\R^n$ and $\bS^n$ and relate them via the conformal
compactification $\R^n \into \bS^n$. This connects the kernel 
$Q(x,y)=(1-\la x , y\ra)^{-s/2}$ on $\bS^n$ to the well known kernel $\|x-y\|^{-s}$ on $\R^n$. 
We show that, for $s = 0$ and $\max(0,n-2) \leq s < n$, this leads to a reflection positive 
distribution cyclic representation $(\pi_s,\cE_s, \alpha)$. Here $\alpha$ can be represented by a 
point measure $\delta_x$ on the  equator corresponding to the unit sphere $\bS^{n-1}$ in 
$\R^n$, the involution $\tau$ corresponds to the reflection in this sphere $\bS^{n-1}$ and 
$S$ is the compression semigroup of the unit ball. 
A Cayley transform translates this into the time reflection 
$(x_0,\ldots ,x_{n-1})\mapsto (-x_0,x_1,\ldots ,x_{n-1})$ and transforms the open unit ball into the 
open half space and the semigroup into the conformal compression semigroup of the half space. 
In \cite[Lemma~2.1]{FL10} and \cite[Lemma~3.1]{FL11}, 
Frank and Lieb give two different proofs of the reflection positivity 
of the distribution $\|x\|^{-s}$, $\max(0,n-2) \leq s < n$ on $\R^n$ 
with respect to reflections in a half space. The also use conformal invariance 
of the corresponding kernel to obtain results similar to our Theorem~\ref{thm:6.7}.


It is worth pointing out that we have a tower of groups
\[
\begin{matrix}
\OO^+_{1,n+1}&\longleftrightarrow & \OO_{2,n}\\
\bigcup & & \bigcup\\
\cE_{n} & \longleftrightarrow &\cP_{n}\\
\bigcup & & \bigcup\\
\R^{n}&\longleftrightarrow & \R^{1, n-1}
\end{matrix}
\]
where the horizontal arrows stand for $c$-duality and the vertical $\bigcup$ stands 
for inclusions. Hence the reflection
positivity on the top line results in reflection positivity on each of the other levels.  Reflection positivity at the top
level is quite rare as our restriction of the parameter,  $s=0$ and
$n-2\le s<n$, shows. In particular, our construction does not allow for reflection positivity on 
the direct limit group $\OO^+_{1,\infty}$. On
the other hand, reflection positivity on the bottom 
line is quite common. Much less is known about the physically interesting
part in the middle (cf.\ \cite{KL82}, \cite{Jo86, Jo87}). 

We would also like to point out that neither we 
nor the previous articles \cite{JO98,JO00} discuss the interesting case of vector-valued
complementary series representations. 
Our results on vector-valued kernels should provide some of the techniques to treat this case. 

{\bf Acknowledgement:} We thank R.~Frank and E.~Lieb for pointing 
out the references \cite{FL10,FL11} and \cite{FILS78}. 

\tableofcontents

\subsection*{Notation}

We write $\R_+ := ]0,\infty[$ for the positive open half line.
Elements of $\R^{n}$, $n \in \N$, are written
$x = (x_0,x_1, \ldots, x_{n-1})$, and
$\R^n_+ := \{ x \in \R^n \: x_0 > 0\}$ is the open half space.
The euclidean inner
product on $\R^{n}$ is denoted $\la x,y \ra = \sum_{j = 0}^{n-1}
x_j y_j$, and
\begin{equation}
  \label{eq:lormet}
[x,y] := x_0 y_0 - x_1 y_1 - \cdots -x_{n-1} y_{n-1}
\end{equation}
is the canonical Lorentzian form on $\R^{n}$, turning it into
the $n$-dimensional Minkowski space.

For a function $\phi \: G \to \C$ on a Lie group $G$, we use the notation
\[ \tilde\phi(g) := \oline{\phi(g^{-1})} \quad \mbox{ and } \quad
\phi^*(g) := \oline{\phi(g^{-1})}\Delta_G(g)^{-1},\]
where $\Delta_G$ is the modular function of $G$, defined by 
\begin{equation}\label{eq:modfun}
 \Delta_G(y) \int_G f(xy)\, d\mu_G(x) =  \int_G f(x)\, d\mu_G(x)
\quad \mbox{ for } \quad f \in C_c(G), y \in G.
\end{equation}
If $\tau$ is an involutive automorphism of $G$, then we also put
\[ g^\sharp := \tau(g)^{-1} \quad \mbox{ and } \quad
\phi^\sharp := \phi^* \circ \tau\, .\]

For the {\it Fourier transform of a measure $\mu$} on the
dual $V^*$ of a finite-dimensional real vector space $V$, we write
\begin{equation}
  \label{eq:ftmeas}
\hat\mu(x) := \int_{V^*}e^{-i\alpha(x)}\, d\mu(\alpha).
\end{equation}
The {\it Fourier transform of an $L^1$-function $f$} on $\R^n$ is defined by
\begin{equation}
  \label{eq:ftfunc}
\hat f(\xi) := (2\pi)^{-n/2} \int_{\R^n} f(x) e^{-i\la \xi,x\ra}\,dx.
\end{equation}
For an involution $\tau$ on $\R^n$, we 
write elements of the dual space as $(\alpha_+,\alpha_-)$
with $\alpha_\pm \circ \tau = \pm\alpha_\pm$
and define the corresponding {\it Fourier--Laplace transform} of a 
 measure $\mu$ on $\R^n$ by
 \begin{equation}
   \label{eq:ftlapl}
{\cF\cL}(\mu)(x) := \int_{\R^n} e^{-i\alpha_+(x)} e^{-\alpha_-(x)}\,
d\mu(\alpha_+, \alpha_-).
 \end{equation}
For $\tau = \id$ this is the Fourier transform $\cF(\mu)=\hat\mu$, and
for $\tau = -\id$, this is the {\it Laplace transform} $\cL(\mu)$.

If $M$ is a smooth manifold and $\cD(M) = C^\infty_c(M,\C)$ is the space
of smooth compactly supported functions on $M$, endowed with its natural
LF topology (\cite{Tr67}), then we write
$\cD'(M)$ for the space of continuous {\it antilinear} functionals on $\cD(M)$;
the {\it distributions} on $M$. We endow this space with the {\it strong dual topology}, i.e., 
the topology of uniform convergence on bounded subsets of $\cD(M)$.

\section{An abstract approach to reflection positivity}\mlabel{sec:1}

In this section we introduce a general concept of a reflection positive 
representation for a triple $(G, \tau,S)$ consisting of a group $G$, 
an involution $\tau$ on $G$ and a subsemigroup $S \subeq G$. 
To study reflection positive representations, it is crucial to control 
the way they are generated by simpler data. This leads to the concept 
of a reflection positive function. 

\subsection{Reflection positive representations} 

A {\it symmetric group} is a pair $(G,\tau)$, consisting of a group $G$
and an involutive automorphism $\tau$ of $G$, which is also allowed
to be trivial, i.e., $\tau = \id_G$. Let $(G,\tau)$ be a symmetric group and $S \subeq G$ be
a subsemigroup which is invariant
under the involution $g \mapsto g^\sharp :=\tau(g)^{-1}$. Then $(S,\sharp)$
is an involutive semigroup, i.e.,
\[ (s^\sharp)^\sharp = s \quad \mbox{ and } \quad
(st)^\sharp = t^\sharp s^\sharp \quad \mbox{ for } \quad s,t \in S.\]
In the following we will write $G_\tau := G \rtimes \{\1,\tau\}$.

\begin{ex}
  \mlabel{ex:1.1}
(a) Let $(G,\tau)$ be a symmetric finite-dimensional Lie group,
$H$ an open subgroup of $G^\tau := \{ g \in G \: \tau(g) = g\}$ and 
$\g = \fh \oplus \fq$ the decomposition of $\L(G) = \g$ into $\L(\tau)$-eigenspaces,
\[ \fh = \g^\tau = \ker(\L(\tau) - \1), \quad \fq = \g^{-\tau} =
\ker(\L(\tau) + \1). \]
If $S \subeq G$ is a subsemigroup of the form
$S = H \exp(C)$, where $C \subeq \fq$ is an $\Ad(H)$-invariant convex cone
(cf.~\cite{HO96}), then the invariance of $S$ under $\sharp$ follows from
$(h\exp x)^\sharp = \exp x h^{-1} = h^{-1}\exp(\Ad(h)x)$ for
$h \in H$, $x \in \fq$.

(b) If $G$ is an abelian group and $\tau(g) = g^{-1}$, then
every subsemigroup of $G$ is invariant under the involution
$\sharp = \id$.

Below we shall also encounter finite-dimensional vector spaces
$G = (V,+)$ with an involution $\tau$ and open convex cones
$\Omega \subeq V$ invariant under $\sharp$.
\end{ex}

\begin{defn} Let $\cE$ be a Hilbert space and 
$\theta \in \U(\cE)$ an involution. 
We call a closed subspace $\cE_+ \subeq \cE$ {\it $\theta$-positive} 
if $\la v,v\ra_\theta := \la \theta v,v \ra \geq 0$ for $v \in \cE_+$. 
We then say that  the triple $(\cE,\cE_+,\theta)$ is a {\it reflection positive 
Hilbert space}. In this case we write 
\[ \cN := \{ v \in \cE_+ \: \la \theta v, v \ra = 0\},\] 
$q \: \cE_+ \to \cE_+/\cN, v \mapsto [v]= q(v)$ for the quotient map 
and $\hat\cE$ for the Hilbert completion of $\cE_+/\cN$ with respect to 
the norm $\|[v]\| := \sqrt{\la \theta v, v \ra}$. 
\end{defn}

\begin{defn}
  \mlabel{def:1.2}
Let $(\cE,\cE_+, \theta)$ be a reflection positive Hilbert space. 
A unitary representation $(\pi,\cE)$ of $G$ is said to be
{\it reflection positive} with respect to the triple $(G,\tau,S)$ 
if it extends to a unitary representation $\pi$ of $G_\tau$ with 
$\pi(\tau) = \theta$ and $\pi(S)\cE_+ \subeq \cE_+$. 
Note that the extendibility of a unitary representation of $G$ to $G_\tau$ is equivalent to the 
existence of a unitary involution $\theta$ on $\cE$ satisfying 
$\theta \pi(g) \theta = \pi(\tau(g))$ for $g \in G$. 
\end{defn}

\begin{lem}
  \mlabel{lem:1.7}
If  $(\pi, {\cal H})$ is a reflection positive representation of $G$ on 
$(\cE,\cE_+, \theta)$, then $\hat\pi(s)[v] := [\pi(s)v]$ defines a representation
$(\hat\pi, \hat\cE)$ of the involutive semigroup $(S,\sharp)$ by contractions.
\end{lem}

\begin{prf} (cf.\ \cite{JO00})
For $s \in S$ and $v,w \in \cE_+$ we have
\begin{align*}
\la \pi(s)v,w \ra_\theta
&= \la \theta\pi(s)v, w\ra
= \la v, \pi(s^{-1})\theta w\ra
= \la v, \theta\pi(s^\sharp) w\ra
= \la v, \pi(s^\sharp) w\ra_\theta,
\end{align*}
so that the representation of $S$ on $\cE_+$
is involutive with respect to the hermitian form $\la \cdot, \cdot \ra_\theta$.
Further,
\[ \la \pi(s)v,\pi(s)v\ra_\theta = \la \theta\pi(s)v,\pi(s)v\ra  \leq \|v\| \]
holds for every $s \in S$,
so that \cite[Lemma~II.3.8]{Ne00} implies that the action of
$S$ on $\cE_+$ induces a contraction representation $(\hat\pi,\hat\cE)$ of
$(S,\sharp)$ on the Hilbert space $\hat\cE$.
\end{prf}

\subsection{Operator-valued positive definite functions} 
\mlabel{subsec:d.1}

A serious problem of the concept of a reflection positive representation is 
that it does not behave well under
direct product decompositions in the sense that if
a reflection positive representation
$\pi$ decomposes as $\pi_1 \oplus \pi_2$, then $\pi_1$ and $\pi_2$
need not be reflection positive because there may be no subspace $\cE_+$
which  is adapted to this decomposition. However,
every non-zero element $v \in \cE_+$ generates a cyclic
reflection positive representation, and since it is natural to study
only reflection positive representations generated by
$\cE_+$, we first focus on on reflection positive representations 
generated in a rather controlled fashion by a subspace $\cF \subeq \cE_+$ 
(cf.\ Definition~\ref{def:1.2b}). First we recall some facts and basic 
concepts on positive definite kernels and functions. 

\begin{defn} \mlabel{def:1.5} Let $X$ be a set and $\cF$ be a complex Hilbert space.

\par (a)  A function $K \: X \times X \to B(\cF)$ is called a 
{\it $B(\cF)$-valued kernel}. 
It is said to be {\it hermitian} if 
$K(z,w)^* = K(w,z)$ holds for all $z, w \in X$. 

\par (b) A $B(\cF)$-valued kernel $K$ on $X$ is said to be 
{\it positive definite} if, 
for every finite sequence $(x_1, v_1), \ldots, (x_n,v_n)$ in $X \times \cF$, 
\[ \sum_{j,k = 1}^n \la K(x_j, x_k)v_k, v_j \ra \geq 0. \] 

\par (c) If $(S,*)$ is an involutive semigroup, then a 
function $\phi \: S \to B(\cF)$ is called {\it positive definite} 
if the kernel $K_\phi(s,t) := \phi(st^*)$ is positive definite. 
\end{defn}

Positive definite kernels can be characterized as those 
for which there exists a Hilbert space $\cH$ and a 
function $\gamma \: X \to B(\cH,\cF)$ such that 
\begin{equation}
  \label{eq:evalprod}
K(x,y) = \gamma(x)\gamma(y)^* \quad \mbox{ for } \quad x,y \in X 
\end{equation}
(cf.\ \cite[Thm.~I.1.4]{Ne00}). 
Here one may assume that the vectors 
$\gamma(x)^*v$, $x \in X, v \in \cF$, span a dense subspace of 
$\cH$. If this is the case, then the pair $(\gamma,\cH)$ is called a {\it realization 
of $K$}. 
The map $\Phi \: \cH \to \cF^X, \Phi(v)(x) := \gamma(x)v$, 
then realizes $\cH$ as a Hilbert subspace of $\cF^X$ 
with continuous point evaluations $\ev_x \: \cH \to \cF, f \mapsto f(x)$. 
Then $\Phi(\cH)$ is the unique Hilbert space in $\cF^X$ with continuous point evaluations 
$\ev_x$, for which $K(x,y) = \ev_x \ev_y^*$ for $x,y \in X$. 
We write $\cH_K \subeq \cF^X$ for this subspace and call it 
the {\it reproducing kernel Hilbert space with kernel~$K$}.

\begin{ex} \mlabel{ex:vv-gns} (Vector-valued GNS construction) 
(cf.\ \cite[Sect.~3.1]{Ne00}) Let $(\pi, \cH)$ be a representation of the 
unital involutive semigroup $(S,*)$, $\cF \subeq \cH$ be a closed subspace for which 
$\pi(S)\cF$ is total in $\cH$ and $P \: \cH \to \cF$ denote the orthogonal projection. 
Then $\phi(s) := P\pi(s)P^*$ is a $B(\cF)$-valued positive definite function 
on $S$  with $\phi(\1) = \1_\cF$ because $\gamma(s) := P\pi(s) \in B(\cH,\cF)$ 
satisfies 
\[ \gamma(s)\gamma(t)^* = P\pi(st^*)P^* = \phi(st^*).\] 
The map 
\[ \Phi \: \cH \to \cF^S, \quad \Phi(v)(s) = \gamma(s)v = P \pi(s)v \] 
is an $S$-equivariant realization of $\cH$ as the reproducing kernel space
$\cH_\phi \subeq \cF^S$, on which $S$ acts by right translation, i.e., 
$(\pi_\phi(s)f)(t) = f(ts)$. 

Conversely, let $S$ be a unital involutive semigroup 
and $\phi \: S \to B(\cF)$ be a positive definite function with 
$\phi(\1) = \1_\cF$. 
Write $\cH_\phi \subeq \cF^S$ for the corresponding reproducing kernel space and 
$\cH_\phi^0$ for the dense subspace spanned by 
$\ev_s^*v, s \in S, v \in \cF$. 
Then $(\pi_\phi(s)f)(t) := f(ts)$ defines a 
$*$-representation of $S$ on $\cH_\phi^0$. 
We call $\phi$ is {\it exponentially bounded} if 
all  operators $\pi_\phi(s)$ are bounded, so that we actually 
obtain a representation of $S$ by bounded operators on $\cH_\phi$. 
As $\1_\cF = \phi(\1) = \ev_\1 \ev_\1^*$, the map 
$\ev_\1^* \: \cF \to \cH$ is an isometric inclusion, so that we may identify 
$\cF$ with a subspace of $\cH$. Then $\ev_\1 \: \cH \to \cF$ corresponds to the 
orthogonal projection onto $\cF$ and 
$\ev_\1 \circ \pi_\phi(s) = \ev_s$ leads to 
\begin{equation}
  \label{eq:factori}
\phi(s) = \ev_s \ev_\1^* =\ev_\1 \pi_\phi(s) \ev_\1^*.
\end{equation}

If $S = G$ is a group with $s^* = s^{-1}$, then $\phi$ is always exponentially bounded and the 
representation $(\pi_\phi, \cH_\phi)$ is unitary. 
\end{ex}

\begin{ex} \mlabel{ex:ker-refpos}
(a) Let $(G,\tau,S)$ be as above and 
$X$ a set on which $G_\tau$ acts. 
We assume that $K \: X \times X \to \C$ is a positive definite 
kernel for which $G_\tau$ acts unitarily on the corresponding 
reproducing kernel space $\cE := \cH_K \subeq \C^X$ by 
\[ (\pi(g)f)(x) = J_g(x)f(g^{-1}.x),\]  
where $J_g \: X \to \C^\times$ is a function with $J_\tau = 1$. For the functions 
$K_x(y) := K(y,x)$ in $\cE$ we then have 
\begin{equation}
  \label{eq:kernshift}
\pi(g) K_x = \oline{J_{g^{-1}}(x)} K_{g.x} \quad \mbox{ for } \quad g \in G, x \in X 
\end{equation}
(cf.\ \cite[Lemma~II.4.1]{Ne00}).

Suppose that $\cD \subeq X$ is an $S$-invariant subset with the 
property that the kernel 
$K^\tau(x,y) := K(x,\tau y)$ is positive definite on $\cD$. 
Then we consider the closed subspace 
\[ \cE_+ := \oline{\Spann \{ K_x \: x \in \cD\}} \] 
and observe that $\cE_+$ is invariant under $S$ by \eqref{eq:kernshift}. From the relation 
\[ \la \pi(\tau)K_y, K_x \ra 
= \la K_{\tau y}, K_x \ra = K(x,\tau y) = K^\tau(x,y)\] 
it now follows that $\cE_+$ is a $\theta$-positive subspace for the 
involution $\theta := \pi(\tau)$. 

(b) In some cases the semigroup $S$ can be obtained from $\cD$. 
Suppose that $X$ is a topological space, that $G_\tau$ acts continuously 
and that $\cD \subeq X$ is open with $\tau(\cD) = (X\setminus \cD)^0$. 
Then the semigroup 
\[ S_\cD := \{ g \in G \: g.\cD \subeq \cD\} \] 
is invariant under $\sharp$. In fact, $g\cD \subeq \cD$ implies that 
$\tau(g)\tau\cD \subeq \tau\cD$ and hence 
$g^\sharp (\tau\cD) \supeq \tau\cD$. But as $\cD = (X \setminus \tau\cD)^0$, 
this leads to $g^\sharp.\cD \subeq \cD$. 
\end{ex}

\subsection{Reflection positive functions}

After these preparations, we now turn to a suitable 
concept of positive definite functions compatible with reflection 
positivity.

\begin{defn}\mlabel{def:1.2e}  
Let $(G,\tau)$ be a symmetric group, $S \subeq G$ a subsemigroup 
invariant under $s^\sharp := \tau(s)^{-1}$ and 
$\cF$ be a Hilbert space. We call a function 
$\phi: G \to B(\cF)$ {\it reflection positive} (with respect to $S$) 
if the following conditions are satisfied: 
\begin{description}
\item[\rm(RP1)] $\phi$ is positive definite, 
\item[\rm(RP2)] $\phi\circ \tau = \phi$, and 
\item[\rm(RP3)] $\phi\res_S$ is positive definite 
as a function on the involutive semigroup $(S,\sharp)$.
\end{description}
\end{defn} 

\begin{defn}
  \mlabel{def:1.2b}  
A triple $(\pi,\cE,\cF)$, where $(\pi,\cE)$ is a 
unitary representation of $G_\tau$ and $\cF \subeq \cE$ 
is a $G$-cyclic subspace fixed pointwise by $\theta := \pi(\tau)$ 
is said to be a 
{\it a reflection positive $\cF$-cyclic representation} 
if the closed subspace $\cE_+ := \oline{\Spann \pi(S)\cF}$ is 
$\theta$-positive. 

If, in addition, $\cF = \C v_o$ is one-dimensional, 
then we call the triple $(\pi,\cE,v_o)$ 
{\it a reflection positive cyclic representation}. 
\end{defn}

\begin{rem} Suppose that $G$ is a topological group and 
$S \subeq G$ is such that $\1 \in \oline S$. 
Let $(\pi,\cE)$ be a continuous unitary representation of $G_\tau$,  $\cF \subeq \cE$ 
be a closed $\theta$-invariant subspace and 
$\cE_+ := \oline{\Spann \pi(S)\cF}$. 
Then $\pi(S)\cF \subeq \cE_+$ implies that $\cF \subeq \cE_+$. 
As $\cF$ is $\theta$-invariant, each element $v \in \cF$ can be written 
as $v = v_+ + v_-$ with $\theta v = v_+ - v_-$. 
If $\cE_+$ is $\theta$-positive, then $\cF$ is also $\theta$-positive and we obtain 
$\theta v = v$ for each $v \in \cF$. Therefore, in this context, the assumption that $\theta$ fixes every element 
of $\cF$ is equivalent to the $\theta$-invariance of $\cF$. 
\end{rem}

\begin{prop}
  \mlabel{prop:1.xx} {\rm(Reflection positive GNS construction)} 
Let $(G,\tau)$ be a symmetric group and $S \subeq G$ be a $\sharp$-invariant subsemigroup. 
\begin{description}
\item[\rm(i)] If $(\pi, \cE,\cF)$ is an $\cF$-cyclic reflection positive 
representation of $G_\tau$ and $P \: \cE \to \cF$ the orthogonal projection, 
then $\phi(g) := P \pi(g) P^*$ is a reflection positive function 
on~$G$ with $\phi(\1) = \1_\cF$. 
\item[\rm(ii)] Let $\phi \: G \to B(\cF)$ is a reflection positive function on $G$ 
with $\phi(\1) = \1_\cF$ and $\cH_\phi \subeq \cF^G$ be  the Hilbert subspace 
with reproducing kernel $K(x,y) := \phi(xy^{-1})$ on which $G$ acts by $(\pi_\phi(g)f)(x) := f(xg)$ 
and $\tau$ by $(\tau f)(x) := f(\tau(x))$. 
We identify $\cF$ with the subspace $\ev_\1^*\cF \subeq \cH_\phi$. 
Then $(\pi_\phi, \cH_\phi,\cF)$ is an $\cF$-cyclic reflection positive 
representation and we have an $S$-equivariant unitary map 
\[ \Gamma \: \hat\cE \to \cH_{\phi\res_S}, \quad 
\Gamma([f]) = f\res_S.\] 
\end{description} 
\end{prop}

\begin{prf} (i) Clearly, $\phi$ is positive definite with $\phi(\1) = \1$. 
That $\phi$ is $\tau$-invariant follows from $\theta\res_\cF = \id_\cF$, which leads to 
$\theta P^* = P^*$ and thus to $P = P\theta$: 
\[ \phi(\tau g \tau) 
= P \theta \pi(g) \theta P^* = P \pi(g) P = \phi(g).\] 
For $s \in S$ we have 
\[  \phi(st^\sharp) 
= P \pi(s) \theta\pi(t^{-1}) \theta P^*  
= P \pi(s) \theta \cdot \theta \pi(t)^* P^*, \] 
so that 
\[ \phi(st^\sharp) = \gamma(s) \gamma(t)^* \quad \mbox{ for } \quad 
\gamma(s) = P \pi(s)\theta.\] 
Therefore $\phi$ is positive definite. 

(ii) For the evaluation maps $\ev_x \: \cH_\phi \to \cF, f \mapsto f(x)$ we have 
$\ev_x \pi_\phi(g)  = \ev_{xg}$. Therefore 
$\pi_\phi(g) \cF = \pi_\phi(g) \ev_\1^* \cF = \ev_{g^{-1}}^*\cF$ shows that 
$\cF$ is $G$-cyclic in $\cH_\phi$. 

For $v \in \cF$, the corresponding element $f \in \cH_\phi$ is the function 
$f(g) = \ev_g \ev_\1^*v= \phi(g)v$. Then 
\[ (\theta f)(g) = f(\tau(g)) = \ev_{\tau(g)} \ev_\1^* v = \phi(\tau(g))v = \phi(g) v 
= f(g) \] 
shows that $\theta\res_{\cF} = \id_\cF$. 

To see that $\cE_+ := \oline{\Spann \pi_\phi(S)\cF}$ is $\theta$-positive, we note that, for $v,w \in \cF$ 
and $s, t\in S$, we have  
\begin{align*}
& \la \theta \pi_\phi(s)\ev_\1^* v, \pi_\phi(t) \ev_\1^* w \ra 
=  \la \pi_\phi(t)^* \theta \pi_\phi(s)\ev_\1^* v, \ev_\1^* w \ra 
=  \la  \pi_\phi(t^\sharp s) \ev_\1^* v, \theta \ev_\1^* w \ra \\
&=  \la  \ev_{s^{-1}\tau(t)}^* v,  \ev_\1^* w \ra 
=  \la  \ev_\1\ev_{s^{-1}\tau(t)}^* v, w \ra 
= \la \phi(t^\sharp s)v,w\ra.
\end{align*}
As $\phi\res_S$ is positive definite on $(S,\sharp)$, it follows that 
$\cE_+$ is $\theta$-positive.

Since $\theta\res_\cF = \id_\cF$, the hermitian form $\la \cdot, \cdot \ra_\theta$ 
restricts on $\cF$ to the 
original scalar product, so that we obtain an isometric inclusion 
$\iota \: \cF \to \hat\cE$. As 
$\hat\pi_\phi(S) \iota(\cF) = [\pi_\phi(S)\cF]$, the subspace $\iota(\cF)$ is $S$-cyclic in $\hat\cE$, 
and for $s \in S$ and $v,w \in \cF$ we have 
\[ \la \iota^* \hat\pi_\phi(s) \iota(v), w \ra 
= \la \hat\pi_\phi(s)\iota(v), \iota(w) \ra = \la \theta\pi_\phi(s)v, w \ra 
= \la \pi_\phi(s)v, w\ra,\] 
so that 
\[ \iota^* \hat\pi_\phi(s) \iota = \ev_\1 \pi_\phi(s) \ev_\1^* 
= \ev_s \ev_\1^* = \phi(s).\] 
This proves that the map 
\[ \Gamma \: \hat\cE \to \cF^S, \quad \Gamma(v)(s) := \iota^* \hat\pi_\phi(s) v \] 
defines an $S$-equivariant isomorphism $\hat\cE \to \cH_{\phi\res_S}.$ 

For $v \in \cF$ and $s \in S$ we further 
have 
\[ \Gamma([\pi_\phi(s)v])(t) 
= \iota^* \hat\pi_\phi(t)[\pi_\phi(s)v] 
= \iota^* [\pi_\phi(ts)v] = \ev_\1 \pi_\phi(ts)v = \ev_t \pi_\phi(s)v,\] 
which implies that 
\[ \Gamma([f])(t) = f(t) \quad \mbox{ for } \quad f \in \cE_+, t \in S.\] 
Hence the natural map $\cH_\phi \supeq \cE_+ \to \cH_{\phi\res_S}$ is simply given by 
restriction to $S$. 
\end{prf}

\begin{cor}
  \mlabel{cor:1.x}
Let $(G,\tau)$ be a symmetric group and $S \subeq G$ be a $\sharp$-invariant subsemigroup. 
\begin{description}
\item[\rm(i)] If $(\pi, \cE,v)$ is a cyclic reflection positive
representation of $G_\tau$,
then $\pi^v(g) := \la \pi(g)v,v\ra$ is a reflection positive function
on~$G$.
\item[\rm(ii)] If $\phi$ is a reflection positive function on $G$,
then $(\pi_\phi, \cH_\phi,\phi)$ is a cyclic reflection positive
representation. 
\end{description}
\end{cor}

\section{Reflection positive distributions} \mlabel{sec:2}

As we shall see below, in some cases we have to pass from
reflection positive functions to the more general class of
reflection positive distributions. This applies
in particular to representations that are most naturally
realized in spaces  of distributions on
a homogeneous space $G/H$. For non-compact subgroups $H$,
we are thus forced to consider unitary representations with $H$-invariant
distribution vectors. Accordingly, we have to study reflection
positive distributions in addition to reflection positive
functions.

\subsection{Positive definite distributions}

In the following we write $\cD(M) = C^\infty_c(M,\C)$ for the space
of compactly supported smooth functions of a manifold $M$ and endow
this space with the usual LF topology, i.e., the locally convex direct limit
of the Fr\'echet spaces $\cD_X(M)$ of test functions supported in
the compact subset $X \subeq M$ (cf.\ \cite{Tr67}). Its antidual, i.e., the
space of continuous antilinear functionals on $\cD(M)$ is the
space $\cD'(M)$ of {\it distributions on $M$}.

\begin{defn} \mlabel{def:2.1}
(a) If $G$ is a Lie group, then $\cD(G)$ is an involutive
algebra with respect to the convolution product and
$\phi^*(g) := \oline{\phi(g^{-1})} \Delta_G(g^{-1})$, 
where $\Delta_G$ is the modular function from \eqref{eq:modfun}. 
Accordingly, we call a distribution $D \in \cD'(G)$  {\it positive
definite}, if it is a positive functional on this algebra, i.e.,
\begin{equation}
  \label{eq:distposdef}
D(\phi^* * \phi) \geq 0 \quad \mbox{ for } \quad \phi \in \cD(G).
\end{equation}

(b) If $\tau$ is an involution on $G$ and $S \subeq G$ an open subsemigroup
invariant under $s \mapsto s^\sharp := \tau(s)^{-1}$, then
$\cD(S)$ is a $*$-algebra with respect to the convolution product
and the $*$-operation $\phi^\sharp := \phi^* \circ \tau$. Accordingly,
we call a distribution $D \in \cD'(S)$ {\it positive definite} if
\begin{equation}
  \label{eq:distposdef2}
D(\phi^\sharp * \phi) \geq 0 \quad \mbox{ for } \quad
\phi \in \cD(S).
\end{equation}
\end{defn}

\begin{rem} \mlabel{rem:2.2}
If $D_h(\phi) = \int_G \oline{\phi(g)}h(g)\, d\mu_G(g)$ holds for a 
locally integrable function $h$ on $G$, then 
\begin{align*}
D_h(\phi^* * \phi) 
&= \int_G \int_G {\phi(x^{-1})}\Delta_G(x)^{-1} 
\oline{\phi(x^{-1}y)} h(y)\, d\mu_G(x)\, d\mu_G(y)\\
&= \int_G \int_G {\phi(x)}\oline{\phi(xy)} h(y)
\, d\mu_G(x) \, d\mu_G(y)\\
&= \int_G \int_G {\phi(x)}\oline{\phi(y)} h(x^{-1}y)
\, d\mu_G(x) \, d\mu_G(y).   
\end{align*}

If $h$ is continuous and positive definite, this formula implies that 
$D_h$ is a positive definite distribution. One can easily see 
that, conversely, $h$ is positive definite if $D_h$~is. 
\end{rem}

\begin{defn} (Distribution vectors)
Let $(\pi, {\cal H})$ be a continuous unitary
representation of the Lie group $G$ on the
Hilbert space ${\cal H}$. We write $\cH^\infty$ for the linear subspace of 
smooth vectors, i.e., of all elements $v \in \cH$ for which the orbit map
$\pi^v \: G \to \cH, g \mapsto \pi(g)v$ is smooth. Identifying
$\cH^\infty$ with the closed subspace of equivariant maps
in the Fr\'echet space $C^\infty(G,\cH)$, we obtain a natural Fr\'echet
space structure on $\cH^\infty$ for which the
$G$-action on this space is smooth
and the inclusion $\cH^\infty \to \cH$ (corresponding to evaluation in $\1 \in
G$) is a continuous linear map  (cf.\ \cite{Mag92}, \cite{Ne10}).

We write ${\cal H}^{-\infty}$ for the space of continuous antilinear
functionals on ${\cal H}^\infty$, the space of {\it distribution vectors},
and note that we have a natural linear
embedding $\cH \into \cH^{-\infty}, v \mapsto \la v, \cdot \ra$.
Accordingly, we also write $\la \alpha, v \ra 
= \oline{\la v, \alpha\ra}$ for $\alpha(v)$ for 
$\alpha \in \cH^{-\infty}$ and $v \in \cH^\infty$.
The group $G$ acts naturally on $\cH^{-\infty}$ by
\[ (\pi^{-\infty}(g)\alpha)(v) := \alpha(\pi(g)^{-1}v),\]
so that we obtain a $G$-equivariant chain of continuous inclusions
\begin{equation}\label{eq:rig}
  {\cal H}^\infty \subeq {\cal H} \subeq {\cal H}^{-\infty} 
\end{equation}
(cf.\ \cite[Sect.~8.2]{vD09}). It is $\cD(G)$-equivariant,
if we define the representation of $\cD(G)$ on $\cH^{-\infty}$ by
\begin{align*}
(\pi^{-\infty}(\phi)\alpha)(v)
&:= \int_G \phi(g) \alpha(\pi(g)^{-1}v)\, d\mu_G(g)
= \alpha(\pi(\phi^*)v).
\end{align*}
\end{defn}

\begin{rem} An important point is that, for any $\phi \in \cD(G)$
and $\alpha \in \cH^{-\infty}$, we have $\pi^{-\infty}(\phi)\alpha
\in \cH^\infty$ in the sense of \eqref{eq:rig}. 
To see this, we first note that $\pi(\phi)$ defines a continuous linear map $\cH \to \cH^\infty$, so that its adjoint maps
$\cH^{-\infty}$ into $\cH$ and then apply the
Dixmier--Malliavin Theorem \cite[Thm.~3.1]{DM78}, 
asserting that
\[ \cD(G) = \cD(G) * \cD(G).\]
Further,  the map
\[\gamma_\alpha \: \cD(G) \to \cH^\infty, \quad \phi \mapsto \pi^{-\infty}(\phi)\alpha \]
is continuous because it is continuous on each Fr\'echet space
$\cD_X(G)$, $X \subeq G$ compact, on which it follows from the Closed Graph Theorem and the
continuity of the orbit maps $\cD(G) \to \cH$.
Therefore
\begin{equation}
  \label{eq:vectodist}
\pi^\alpha(\phi) := \la \alpha, \pi^{-\infty}(\phi)\alpha \ra 
\end{equation}
defines a distribution on $G$. That this distribution is positive
definite follows from
\begin{equation}
  \label{eq:posdef}
\pi^\alpha(\phi^* * \phi)
= \la \alpha, \pi^{-\infty}(\phi^* * \phi) \alpha \ra
= \la \pi^{-\infty}(\phi)\alpha, \pi^{-\infty}(\phi)\alpha \ra \geq 0.
\end{equation}
\end{rem}

\begin{defn}
We say that $\alpha \in \cH^{-\infty}$ is {\it cyclic} if
$\pi^{-\infty}(\cD(G))\alpha$ is a dense subspace of $\cH$.
\end{defn}
  
\begin{thm} \mlabel{thm:vD} {\rm(\cite[Thm.~8.2.1]{vD09})} 
Let $(\pi, \cH)$ be a continuous unitary representation of the Lie group 
$G$. 
\begin{description}
\item[\rm(i)] For each distribution vector 
$\alpha \in \cH^{-\infty}$, there exists a $G$-equivariant continuous linear map 
\[ \eta_\alpha \: \cH \to \cD'(G), \quad 
\eta_\alpha(v)(\phi) := \la v, \pi^{-\infty}(\phi) \alpha\ra 
= \la \pi(\phi^*)v, \alpha\ra. \] 
which is injective if and only if $\alpha$ is cyclic. 
This establishes a one-to-one correspondence between distribution 
vectors and $G$-equi\-variant continuous linear maps 
$\cH \to \cD'(G)$. 
\item[\rm(ii)] The assignment 
\[ 
\alpha \mapsto \eta_\alpha^*, \quad 
\eta_\alpha^*(\phi) := \pi^{-\infty}(\phi)\alpha \] 
establishes a one-to-one correspondence between distribution 
vectors and $G$-equi\-variant continuous linear maps $\cD(G) \to \cH$.
\item[\rm(iii)] The map $\eta_\alpha$ extends to 
\[ \tilde\eta_\alpha \: \cH^{-\infty} \to \cD'(G), \quad 
\tilde\eta_\alpha(\beta)(\phi) := \la \beta, \pi^{-\infty}(\phi)\alpha\ra.\] 
\end{description}
\end{thm}

We will now discuss how every positive definite distribution $D \in \cD'(G)$ leads to a
unitary representation $(\pi, \cH)$ with cyclic distribution vector 
$\alpha \in \cH^{-\infty}$ which defines an embedding $\cH \into \cD'(G)$.
First we recall some functional analytic facts from \cite{Tr67}.

\begin{rem} \mlabel{rem:2.3}
LF spaces are {\it barreled,}
 i.e., all closed absolutely convex absorbing subsets (the barrels) are $0$-neighborhoods. 
This applies in particular to $\cD(M)$ for any $\sigma$-compact
manifold $M$ (\cite[p.~347]{Tr67}). It follows from
\cite[Prop.~34.4]{Tr67} that $\cD(M)$ is a {\it Montel space},
which means that it is barreled and every bounded closed subset of $\cD(M)$ is compact.
This implies that $\cD(M)$ is reflexive (\cite[p.~376]{Tr67})
and that every weakly (= weak-$*$) convergent sequence in
$\cD'(M)$ converges (\cite[p.~358]{Tr67}).
\end{rem}

For the following proposition we recall that the vector-valued GNS construction 
(Example~\ref{ex:vv-gns}) yields for every positive definite distribution 
$D \in \cD'(G)$ a corresponding Hilbert space $\cH_D$ which is contained in the space 
$\cD(G)^*$ of all antilinear functionals on $\cD(G)$. 

\begin{prop}\mlabel{prop:2.4}
Let $D \in \cD'(G)$ be a positive definite distribution on the 
Lie group  $G$ and $\cH_D$ be the corresponding 
reproducing kernel Hilbert space with kernel 
$K(\phi,\psi) := D(\psi^* * \phi)$ obtained by completing 
$\cD(G) * D$ with respect to the scalar product 
$\la \psi * D, \phi * D \ra = D(\psi^* * \phi)$. Then the 
following assertions hold: 
\begin{description}
\item[\rm(i)] $\cH_D \subeq \cD'(G)$ and 
the inclusion $\gamma_D \: \cH_D \to \cD'(G)$ is continuous. 
\item[\rm(ii)] We have a unitary representation 
$(\pi_D, \cH_D)$ of $G$ by $\pi_D(g)E = g_*E$, where \break 
$(g_*E)(\phi) := E(\phi \circ \lambda_g)$ and the integrated representation 
of $\cD(G)$ on $\cH_D$ is given by $\pi_D(\phi)E = \phi * E$. 
\item[\rm(iii)]  There exists a unique distribution 
vector $\alpha_D \in \cH_D^{-\infty}$ with 
$\alpha_D(\phi *D) = D(\phi)$ and 
$\pi^{-\infty}(\phi)\alpha_D = \phi * D$ for $\phi \in \cD(G)$.  
It satisfies $\pi^{\alpha_D} = D$. 
\item[\rm(iv)] $\gamma_D$ extends to a $\cD(G)$-equivariant injection 
$\cH_D^{-\infty} \into \cD'(G)$ mapping  $\alpha_D$ to $D$. 
\end{description}
\end{prop}

\begin{prf} (i) The relation 
$K(\psi, \phi) = \oline{K(\phi,\psi)}$ implies
$D(\phi^*) = \oline{D(\phi)}$ for $\phi \in \cD(G)$. The
functionals $\psi * D$, $\psi \in \cD(G)$,
span a pre-Hilbert space
$\cH_D^0 \subeq \cD'(G)$ with inner product $\la \psi * D, \phi * D \ra = K(\phi, \psi)$ 
on which we have a $*$-representation of $\cD(G)$, given by 
\[ \pi_D(\phi)E := \phi * E.\] 
We claim that the corresponding inclusion map
$\gamma_D \: \cH_D^0 \into \cD'(G)$ 
is continuous with respect to the pre-Hilbert structure on $\cH_D^0$. 
Since $\cH_D^0$ is metrizable, it suffices to show that $\gamma_D$ is
sequentially continuous, so that Remark~\ref{rem:2.3}
further implies that it suffices
to verify weak continuity. This follows immediately from
\begin{equation}
  \label{eq:wcont}
\gamma_D(\psi * D)(\phi) = (\psi * D)(\phi) = D(\psi^* * \phi)
\end{equation}
because the multiplication on $\cD(G)$ is separately 
continuous.\begin{footnote}{That the convolution product on $\cD(G)$ is jointly 
continuous for every Lie group $G$ with at most finitely many connected 
components has  been shown recently by Birth and Gl\"ockner in \cite{BG11}.}  
\end{footnote}
Since $\cD'(G)$ is complete (it is the strong dual of an LF space;
\cite[p.~344]{Tr67}), $\gamma_D$ extends to a continuous
linear map $\gamma_D \: \cH_D \to \cD'(G)$ on the Hilbert space completion
$\cH_D$ of $\cH_D^0$. 

(ii) In view of 
\[ (\phi * T)(\psi) = T(\phi^* * \psi)
= \int_G \phi(g) T(\psi \circ \lambda_g)\, d\mu_G(g), \]
the representation of $\cD(G)$ on $\cH_D^0$ corresponds to
the unitary $G$-representation defined by
\[ (\pi_D(g)T)(\psi) := T(\psi \circ \lambda_g),\]
which is the dual of the left regular representation of $G$
on $\cD(G)$. From the unitarity of the $G$-representation 
on $\cH_D^0$ it follows that it extends to a unitary 
representation $(\pi_D, \cH_D)$ whose continuity follows from 
the continuity of the $G$-orbit maps in $\cD(G)$.

(iii), (iv) Clearly,  the inclusion $\gamma_D$ is $G$-equivariant, so that 
Theorem~\ref{thm:vD} implies the existence of a uniquely determined 
cyclic distribution vector $\alpha_D\in \cH_D^{-\infty}$ satisfying 
\[ \gamma_D(E)(\phi) = \la E, \pi^{-\infty}(\phi)\alpha_D \ra 
= \la \phi^* * E, \alpha_D \ra.\] 
For $E = \psi * D$, this leads to 
\[ D(\psi^* * \phi) {\buildrel \eqref{eq:wcont} \over =} \gamma_D(\psi * D)(\phi)  
= \la \phi^* * \psi * D, \alpha_D \ra.\] 
In view of the Dixmier--Malliavin Theorem (\cite[Thm.~3.1]{DM78}), 
we further get 
\begin{equation}\label{eq:betarel2}
\alpha_D(\phi * D) = \oline{D(\phi^*)} = D(\phi).
\end{equation}
We also obtain 
\[ \la \psi * D, \phi * D \ra 
= D(\psi^* * \phi) 
= \la \phi^* * \psi * D, \alpha_D \ra 
=  \la \psi * D, \pi^{-\infty}(\phi)\alpha_D \ra,  \] 
which leads to 
\begin{equation}
  \label{eq:reprel}
\pi^{-\infty}(\phi)\alpha_D = \phi * D.
\end{equation}

The map $\beta_D \: \cD(G) \to \cH_D, \phi \mapsto \phi * D$
is $\cD(G)$-equivariant and continuous\begin{footnote}
{It suffices to verify continuity on the Fr\'echet subspaces $\cD_X(G)$,
where $X \subeq G$ is a compact subset. On these subspaces, the continuity
follows from the Closed Graph Theorem and the continuity of the maps
$\phi \mapsto \la \beta_D(\phi), \psi * D\ra = \la \phi * D, \psi * D\ra= D(\phi^**\psi)$.}
\end{footnote}
and \eqref{eq:betarel2} means that 
\begin{equation}\label{eq:betarel}
\alpha_D \circ  \beta_D  = D.
\end{equation}
Actually $\beta_D$ defines a continuous linear map 
$\cD(G) \to \cH_D^\infty$ with dense range whose adjoint yields an inclusion 
\[ \beta_D^* \: \cH_D^{-\infty} \into \cD'(G),\] 
i.e., the inclusion $\gamma_D \: \cH_D \into \cD'(G)$ even extends to the larger 
space of distribution vectors. Here \eqref{eq:betarel} means that 
$\beta_D^*\alpha_D = D$. Finally, \eqref{eq:betarel2} implies 
\[ \pi^{\alpha_D}(\phi) = \la \alpha_D, \pi^{-\infty}(\phi) \alpha_D \ra 
= \la \alpha_D, \phi * D \ra 
{\buildrel \eqref{eq:betarel2} \over =}  D(\phi),\] 
i.e., $\pi^{\alpha_D} = D$. 
\end{prf}

\begin{ex} \mlabel{ex:2.7}
If $G = V$ is a real vector group (isomorphic to $\R^n$),
then the Bochner--Schwartz Theorem (\cite[Thm.~XVIII, \S VII.9]{Schw73})
asserts that a distribution $D \in \cD'(V)$ is positive definite if and only if
there exists a tempered positive measure $\mu$ on $V^*$ with
$D = \hat\mu$ (the Fourier transform), i.e.,
\[ D(\phi) = \int_{V^*} \oline{\hat\phi}(\alpha)\, d\mu(\alpha)
\quad \mbox{ and } \quad
 D(\psi^* * \phi) = \la \hat\psi, \hat\phi\ra_{L^2(V^*,\mu)}.\]
{}From this is follows that
$\hat\psi \mapsto \psi* D$ extends to a unitary map
$L^2(V^*,\mu) \to \cH_D$. 
If intertwines the unitary representation of $V$ on $L^2(V^*,\mu)$ by 
\[ (U_v f)(\alpha) := e^{-i\alpha(v)} f(\alpha), \quad f \in L^2(V^*,\mu), v \in V, \alpha \in V^*\]  
with the translation action of $V$ on $\cD'(V)$. 

Since $D = \hat \mu \in \cS'(V)$ is a tempered distribution, the construction
in Proposition~\ref{prop:2.4} actually leads to a continuous inclusion
$\gamma_D \:  \cH_D \into \cS'(V)$.
\end{ex}

\subsection{Reflection positivity for distributions}

After this brief discussion of representations on Hilbert subspaces of
$\cD'(G)$,we now turn to reflection positive distributions on Lie groups and
the corresponding unitary representations.

\begin{defn}
  \mlabel{def:3.5}
If $(G,\tau)$ is a symmetric Lie group and $S$ is open and
$\sharp$-invariant, then we call a distribution
 $D \in {\cal D}'(G)$  {\it reflection positive} for 
$(G,\tau,S)$  if the following conditions are satisfied:
\begin{description}
\item[\rm(RP1)] $D$ is positive definite, i.e.,
$D(\phi^* * \phi) \geq 0$ for $\phi \in \cD(G)$.
\item[\rm(RP2)] $\tau D = D$, i.e., $D(\phi \circ \tau) = D(\phi)$ for
$\phi \in \cD(G)$, and
\item[\rm(RP3)] $D\res_S$ is positive definite
as a distribution on the involutive semigroup $(S,\sharp)$, i.e.,
$D(\phi^\sharp * \phi) \geq 0$ for $\phi \in \cD(S)$.
\end{description}
\end{defn}

\begin{defn}
  \mlabel{def:1.2d}
A triple $(\pi,\cH,\alpha)$, where $(\pi,\cH)$ is a
unitary representation of $G_\tau$ and $\alpha_o \in \cH^{-\infty}$
a cyclic distribution vector fixed under $\theta := \pi(\tau)$,
is said to be a
{\it a reflection positive distribution cyclic representation}
if the closed subspace $\cE_+ := \oline{\Spann \pi^{-\infty}(\cD(S))\alpha_o}$ is
$\theta$-positive.
\end{defn}

\begin{prop}
  \mlabel{prop:1.y}
For $(G,\tau,S)$ as above, the following assertions hold:
\begin{description}
\item[\rm(i)] If $(\pi, \cH,\alpha)$ is a distribution cyclic
reflection positive representation of $G_\tau$,
then $\pi^\alpha(\phi) := \alpha(\pi^{-\infty}(\phi)\alpha)$ is a reflection
positive distribution on~$G$.
\item[\rm(ii)] If $D$ is a reflection positive distribution on $G$,
then $(\pi_D, \cH_D,D)$ is a reflection positive  distribution cyclic
representation, where $\tau$ acts on $\cD'(G)$ by
$(\tau E)(\phi) := E(\phi \circ \tau)$.
\end{description}
\end{prop}

\begin{prf} (i) We have already seen in  \eqref{eq:posdef}
that the distribution $\pi^\alpha$ is positive definite.
That $\pi^\alpha$ is $\tau$-invariant follows from the $\theta$-invariance of
$\alpha$:
\begin{align*}
\pi^\alpha(\phi \circ \tau)
&= \la \alpha, \pi^{-\infty}(\phi \circ \tau)\alpha \ra
= \la \alpha, \pi(\tau)\pi^{-\infty}(\phi)\pi^{-\infty}(\tau)\alpha \ra \\
&= \la \alpha, \pi^{-\infty}(\phi)\alpha \ra = \pi^\alpha(\phi).
\end{align*}
For $\phi \in \cD(S)$ we further have
\begin{align*}
 \pi^\alpha(\phi^\sharp * \phi)
&= \la \alpha, \pi^{-\infty}(\phi^\sharp)\pi^{-\infty}(\phi)\alpha \ra
= \la \pi^{-\infty}(\phi \circ \tau)\alpha, \pi^{-\infty}(\phi)\alpha \ra\\
&= \la \theta\pi^{-\infty}(\phi)\theta\alpha, \pi^{-\infty}(\phi)\alpha \ra
= \la \pi^{-\infty}(\phi)\alpha, \pi^{-\infty}(\phi)\alpha\ra_\theta \geq 0,
\end{align*}
so that $\pi^\alpha$ is reflection positive.

(ii) The $\tau$-invariance of $D$ implies the $\tau$-invariance of the
kernel $K$, so that $(\theta E)(\phi) := E(\phi \circ \tau)$ defines a
unitary operator on $\cH_D$ (cf.\ \cite[Rem.~II.4.5]{Ne00}),
and we thus obtain a
unitary representation $\pi_D$ of $G_\tau$ on $\cH_D$ with
$\pi_D(\tau) = \theta $.

In Proposition~\ref{prop:2.4} we already observed that
$D \in \cH_D^{-\infty}$ is a cyclic distribution vector,
and its $\theta $-invariance follows from the $\tau$-invariance of $D$.
Finally, we note that any $\phi \in \cD(S)$ satisfies
\[ \la \theta  K_\phi, K_\phi \ra
= \la K_{\phi \circ \tau}, K_\phi \ra
= K(\phi, \phi \circ \tau) = D(\phi^\sharp * \phi) \geq 0.\]
Therefore the positive definiteness of $D\res_{\cD(S)}$
implies that the closed subspace $\cE_+$ generated by the elements
$K_\phi$, $\phi \in \cD(S)$, is $\theta $-positive.

The existence of $\Gamma$ follows as in Proposition~\ref{prop:1.xx}(ii). 
\end{prf}

\begin{lem} \mlabel{lem:1.3}
For $s \in ]-\infty,n[$,
the function $\|x\|^{-s}$ on $\R^n$ is locally integrable,
hence defines a distribution $D_s \in \cD'(\R^n)$. It is positive definite
if and only if $s \geq 0$.
\end{lem}

\begin{prf} Using polar coordinates, one immediately sees that, for
$r := \|x\|$, the function $r^{-s}$ is locally integrable if and only if
$s < n$.

For, $0 < s < n$, it follows from \cite[Ex.~VII.7.13]{Schw73}
that the Fourier transform of $r^{-s}$ is given by
\[ \cF(r^{-s}) = \pi^{s-n/2}
\frac{\Gamma\big(\frac{n-s}{2}\big)}{\Gamma\big(\frac{s}{2}\big)} r^{s-n}.\]
Since $0 <s < n$ implies $-n < s-n < 0$ and the $\Gamma$-factors are
positive in this range, this formula shows that $\cF(r^{-s})$
is a positive tempered measure, hence that $r^{-s}$ is positive definite.

For $s \leq 0$, the function $r^{-s}$ is continuous and real-valued.
If it is positive definite, it has a maximal value at 
$x = 0$, which is only the case for $s = 0$. Then $r^s = 1$, which
is obviously positive definite.
\end{prf}

\begin{ex}
The preceding lemma  and Proposition~\ref{prop:2.4} imply that
we have a map $\beta : \cD (\R^n) \to \cD^\prime (\R^n)$ given by
\[\beta (\phi)(\psi) 
=\int_{\R^n}\int_{\R^n}\frac{\phi (x)\overline{\psi (y)}}{\|x-y\|^s}\, dxdy.\]
(cf.\ Remark~\ref{rem:2.2}). 
The completion of $\beta(\cD (\R^n))$ is a translation invariant Hilbert subspace 
$\cH_{D_s} \subeq \cD^\prime (\R^n)$ and we even obtain 
an equivariant embedding $\cH_{D_s}^{-\infty} \subeq \cD^\prime (\R^n)$. 

To connect this with our discussion of reflection positivity, let $\tau $ be the ``time'' reflection 
$\tau(x_0,x_1 \ldots ,x_{n-1})=(-x_0,x_1,\ldots ,x_{n-1})$ mapping the open half space
$\R^n_+=\{(x_0,y)\: x_0>0, y\in\R^{n-1}\}$ onto $-\R^n_+$, keeping the 
bounding hyperplane pointwise fixed. The open half-space 
$S=\R^n_+$ is an involutive semigroup with respect to 
the involution $(x_0,y)^{\#}=(x_0,-y)$. 
Then $D_s \in\cH_{D_s}^{-\infty}$ represents a cyclic distribution vector which is 
fixed under $\theta $ because $\|\tau(x)\| = \|x\|$ for $x \in \R^n$. 
We shall see below that 
$D_s$ is reflection positive for 
$s = 0$ and $\max(0,n-2) \leq s < n$ 
(cf.\ Proposition~\ref{prop:wallset}). 
\end{ex}

\section{Reflection positivity on the real line}\mlabel{sec:3}

A particular special case, where it is interesting to study reflection positivity is the
subsemigroup $S = ]0,\infty[ = \R_+ $ of the real line $G =\R$. Here we
consider the involution $\tau(x) =-x$, so that $S$ is invariant
under the involution $s^\sharp = s$. In this case we obtain a complete 
description of the reflection positive functions in terms of 
an integral representation. 

\subsection{Reflection positive operator-valued functions on $\R$}
\mlabel{subsec:d.2}

We start with a characterization of continuous reflection positive 
functions for $(G,\tau,S) = (\R, -\id_\R, \R_+)$. 

\begin{prop}
  \mlabel{prop:2.1b}  Let $\cF$ be a Hilbert space and
$\phi \: \R \to B(\cF)$ be positive definite and strongly continuous. 
Then $\phi$ is reflection positive if and only if
there exists a finite $\Herm(\cF)_+$-valued Borel measure
$Q$ on $[0,\infty[$ such that
\begin{equation}
  \label{eq:2.1b}
  \phi(x) = \int_0^\infty e^{-\lambda |x|}\, dQ(\lambda).
\end{equation}
\end{prop}

\begin{prf} Suppose first that $\phi$ is
reflection positive for $(\R, -\id,\R_+)$ and put $S :=(\R _+,+)$. Then
$\phi_S := \phi\res_S$ is positive definite with respect to the
trivial involution and corresponds to a contraction
representation of $S$ because
$|\la \phi(s)v,v\ra| \leq \la \phi(\1)v,v\ra$
holds for the positive definite functions
$x \mapsto \la \phi(x)v,v\ra$, $v \in \cF$, on $\R$
(\cite[Cor.~III.1.20(ii)]{Ne00}).
Hence there exists a unique finite $\Herm(\cF)_+$-valued Borel measure
$Q$ on $[0,\infty[$ such that
$$ \phi(x) = \int_0^\infty e^{-\lambda x}\, dQ(\lambda)
\quad \mbox{ for } \quad x > 0 $$
(\cite[Cor.~IV.4]{Ne98}).
The Dominated Convergence Theorem and the (strong) continuity of
$\phi$ then imply 
\[ \phi(0) = Q([0,\infty[) = \int_0^\infty \, dQ(\lambda).\]
Now \eqref{eq:2.1b} follows from the fact that $\phi(-x)= \phi(x)^*
= \phi(x)$ holds for $x \geq 0$.

For the converse, we assume that $\phi$ has an integral representation as
in \eqref{eq:2.1b}. This immediately
implies that $\phi\res_S$ is positive definite on $S$
for the involution $s^\sharp = s$ and that
$\phi$ is continuous (\cite[Prop.~II.11]{Ne98}).
We have to show that $\phi$ is positive definite. Since the cone of
positive definite functions is closed under pointwise limits,
it suffices to show that, for $\lambda > 0$
and $A \in \Herm(\cF)_+$, the function
$\phi_\lambda(x) := e^{-\lambda |x|}A$ is positive definite on $\R$.
For $\lambda = 0$ this is trivial, and for $\lambda > 0$ it
follows from
\begin{equation}\label{eq:7}
\phi_\lambda(x) = e^{-\lambda|x|} A
= \int_\R e^{i xy}A\, d\mu_\lambda(y),\quad \mbox{ where } \quad
d\mu_\lambda(y) = {1\over \pi} {\lambda \over \lambda^2 + y^2} dy
\end{equation}
is the Cauchy distribution (\cite[\S 47]{Ba78}).
\end{prf}

\begin{rem} \mlabel{rem:3.2} (a) For the special 
case where $\phi(0) = \1$, strongly continuous reflection positive 
functions are called (OS)-positive covariance functions in \cite{Kl77}, and 
the preceding result specializes to \cite[Rem.~2.7]{Kl77} in the following sense. 

Using \eqref{eq:factori} to write
$\phi(s) = \ev_\1 \circ \pi_\phi(s) \circ \ev_\1^*$
and representing $\pi_\phi$ by a spectral measure $P$ on $[0,\infty[$ as
\[ \pi_\phi(s) = \int_0^\infty e^{-\lambda s}\, dP(\lambda),\]
we obtain the integral representation of $\phi$ with
$Q(E) := \ev_\1 \circ P(E) \circ \ev_\1^*$ for any Borel subset $E \subeq [0,\infty[$. 

(b) The proof above implies in particular that every bounded 
strongly continuous $B(\cF)$-valued positive definite 
function $\phi \: \R_+ \to B(\cF)$ extends by 
$\tilde\phi(x) := \phi(|x|)$ to a positive definite function on 
$\R$. Clearly, $\tilde\phi$ is reflection positive. 
This observation can also be found in \cite[\S I.8.2, p.~30]{SzN70}. 
If $\phi$ is a representation, i.e., a one-parameter semigroup 
of contractions, then the unitary representation 
$\pi_{\tilde\phi}$ of $\R$ on the reproducing kernel Hilbert space 
$\cH_{\tilde\phi}$ is called the {\it minimal unitary dilation} of~$\phi$. 

If $\phi(0) = \1$, then we can identify $\cF$ with a subspace of 
$\cH_{\tilde\phi}$ and $\phi(s) = P \pi_{\tilde\phi}(s) P^*$ holds 
for $s \geq 0$ and the orthogonal projection $P \: \cH_{\tilde\phi} \to \cF$. 
\end{rem}

Specializing the preceding result to $\cF = \C$, we obtain the following integral 
representation. A discrete version for reflection positivity on the group $\Z$ 
can be found in \cite[Prop.~3.2]{FILS78}. 

\begin{cor}
  \mlabel{cor:2.1}  A continuous function $\phi \: \R \to \C$
is reflection positive if and only if it
has an integral representation of the form
\begin{equation}
  \label{eq:2.1}
  \phi(x) = \int_0^\infty e^{-\lambda|x|}\, d\nu(\lambda),
\end{equation}
where $\nu$ is a finite positive Borel measure on $[0,\infty[$.
\end{cor}

\begin{rem} (a) Since
the function $\phi_\lambda(x) = e^{-\lambda |x|}$, $\lambda \geq 0$, is real-valued,
it can also be written as
$$\phi_\lambda(x) = e^{-\lambda|x|}
= \frac{2\lambda}{\pi} \int_0^\infty \cos(xy)\, \frac{dy}{\lambda^2 + y^2}. $$

(b) Using the isomorphism $\exp \: (\R,+) \to (\R _+,\cdot)$,
we also get a description of the reflection positive functions on
$\R _+$ with respect to the involution
$\tau(a) = a^{-1}$ and the subsemigroup $S = ]0,1[$. They are given by
\[
\phi_\lambda(a) =  e^{-\lambda|\log a|}=
\begin{cases}
a^{-\lambda} & \text{for } a \geq 1 \\
a^{\lambda} & \text{for } a \leq 1. \\
\end{cases}\]
\end{rem}

\begin{ex} \mlabel{ex:2.2b}
We take a closer look at the representation associated to the positive definite function
$\phi(x) = e^{-\lambda |x|}$ on $\R$, where $\lambda > 0$.

(a) Let $d\mu(y) = {1\over \pi} {\lambda \over \lambda^2 + y^2} dy$
be the Cauchy distribution from \eqref{eq:7}, 
so that $\phi = \hat\mu$ is the
Fourier transform of $\mu$.
Then we
may also realize $\cH_\phi$ for $\phi(x) = e^{-\lambda|x|}$
as $L^2(\R,\mu)$, where the unitary
isomorphism is given by
$$ \Gamma \: L^2(\R,\mu) \to \cH_\phi, \quad
\Gamma(f)(x) = \la f, \pi(x)1\ra
= \int_\R f(y) e^{-ixy}\, d\mu(y)
= (f\mu)\,\hat{}(x) $$
and $\tau$ acts on $L^2(\R,\mu)$ by $\theta(f)(x) = f(-x)$
(cf.\ Example~\ref{ex:2.7}). Note that $\Gamma(1) = \phi$. 

We consider the closed subspace $\cE_+ \subeq \cH_\phi$
generated by $\pi_\phi(S)\phi$ and note that
$$ \la \theta \pi_\phi(s)\phi, \pi_\phi(t)\phi \ra
=  \la \pi(-t-s)\phi, \phi \ra = \phi(-t-s)
= \phi(t+s) = e^{-\lambda(t+s)}. $$
Hence the corresponding reproducing kernel space on
$S$ is one-dimensional, generated by the character
$e_\lambda$ with $e_\lambda(x) = e^{-\lambda x}$. Moreover, 
$(\pi_\phi, \cH_\phi)$ is the minimal unitary dilation of $e_\lambda$ 
(cf.\ Remark~\ref{rem:3.2}(b)). 

(b) A slightly different realization of $\cH_\phi$, which makes the
scalar product on this space more explicit, is to consider it 
as a completion of the space $C_c(\R,\C)$ of compactly supported
continuous function, endowed with the hermitian form
\begin{align*}
(f,g) \mapsto
&\la \pi_\phi(f)\phi,  \pi_\phi(g)\phi \ra
= \int_\R \int_\R f(x)\oline{g(y)} \la \pi_\phi(x)\phi,  \pi_\phi(y)\phi \ra
\, dx\, dy \\
&= \int_\R \int_\R f(x)\oline{g(y)} \phi(x-y) \, dx \, dy
= \int_\R \int_\R f(x)\oline{g(y)} e^{-\lambda|y-x|} \, dx \, dy
\end{align*}
(cf.\ Remark~\ref{rem:2.2}). 
In this picture, $\cE_+$ corresponds to the subspace generated by
those functions $f \in C_c(\R,\C)$, which are supported by $[0,\infty[$.
For two such functions $f,g$, we obtain
\begin{align*}
&\la \theta f, g\ra
= \int_\R \int_\R f(-x)\oline{g(y)} e^{-\lambda|y-x|} \, dx \, dy
= \int_\R \int_\R f(x)\oline{g(y)} e^{-\lambda|y+x|} \, dx \, dy \\
&= \int_0^\infty \int_0^\infty f(x)\oline{g(y)} e^{-\lambda(y+x)} \, dx \, dy
= \int_0^\infty f(x)e^{-\lambda x}\, dx \cdot
\oline{\int_0^\infty g(y) e^{-\lambda y}\, dy}.
\end{align*}
This formula reflects the fact that the space $\hat\cE$ is one-dimensional
and that the natural map $\cE_+ \to \hat\cE$ can be realized by
$f \mapsto \int_0^\infty f(x)e^{-\lambda x}\, dx$.
\end{ex}

\subsection{Reflection positive distributions on $\R$}

In this subsection we discuss two interesting families of reflection
positive distributions for the triple $(\R, -\id, \R_+)$.
These distributions occur naturally by restriction of the
complementary series representations of $\SL_2(\R)$ (which is locally isomorphic to 
$\OO_{1,2}(\R)_0$) to the subgroups $N$ and $A$ in the Iwasawa decomposition
(cf.\ \cite{JO00}). These representations are discussed in a more 
general context in Sections~\ref{sec:4} and \ref{sec:5} below. 

\begin{ex}
  \mlabel{ex:2.3}
(a) For $0 < s < 1$ the locally integrable function
$|x|^{-s}$ on $\R$ defines a positive definite measure, resp.,
distribution (Lemma~\ref{lem:1.3}). The corresponding Hilbert space~$\cH_s$ 
is the completion of $C_c(\R)$ with respect to the scalar product
\begin{equation}
  \label{eq:inprod1}
\la f, g \ra := \int_\R \int_\R f(x)\oline{g(y)} |x-y|^{-s}\, dxdy
\end{equation}
(cf.~Proposition~\ref{prop:2.4}).

The distribution $|x|^{-s}$ restricts to the function
$x^{-s}$ on $\R_+$ which is positive definite
with respect to the trivial involution because, for any $\alpha > 0$,
\begin{equation}
  \label{eq:lapform}
\frac{1}{\Gamma(\alpha)} \cL(y^{\alpha-1}\, dy)(x) = x^{-\alpha}
\quad \mbox{ for } \quad x > 0.
\end{equation}
We conclude that the distribution
$|x|^{-s}$ is reflection positive for the triple $(\R, -\id_\R,  \R_+)$.

(b) On the Hilbert space $\cH_s$ we also have a unitary representation
of the multiplicative group $A := \R $, given by
\[ (\pi(a)f)(x) := |a|^{\frac{s}{2}-1} f(a^{-1}x)\]
and an involution
$$ (\theta f)(x) := |x|^{s-2} f(x^{-1}) $$
satisfying $\theta \pi(a)\theta  = \pi(a^{-1})$, so that we obtain
for $\tau(a) := a^{-1}$ a representation of $A_\tau$ on $\cH_s$.

Let $\cE_+ \subeq \cH_s$ be the closed subspace generated by functions
supported in $]-1,1[$ and note that it is invariant under the action
of the subsemigroup $S := \{ a \in A \: |a| <  1\}$.
For $f \in C_c(]-1,1[,\C)$, we have
\begin{align*}
\la \theta f,f \ra
&= \int_\R \int_\R |x|^{s-2} f(x^{-1})\oline{f(y)} |x-y|^{-s}\, dx\, dy \\
&= \int_\R \int_\R |x^{-1}|^{s-1} f(x)\oline{f(y)} |x^{-1}-y|^{-s}\, \frac{dx}
{|x|}\, dy \\
&= \int_{-1}^1 \int_{-1}^1 f(x)\oline{f(y)} |x|^{-s} |x^{-1}-y|^{-s}
\, dx\, dy \\
&= \int_{-1}^1 \int_{-1}^1 f(x)\oline{f(y)} (1 -xy)^{-s}
\, dx\, dy \geq 0
\end{align*}
because the kernel $(1-xy)^{-s}= \sum_{n=0}^\infty {n-1+s \choose n} (xy)^n$
on $]-1,1[$ is positive definite, which follows from the
non-negativity of the binomial coefficients.
Therefore the unitary representation $(\pi, \cH_s)$ of $A$ is reflection
positive.

We now describe a cyclic distribution vector and the corresponding
distribution on~$A$.
For the integrated representation we find
\[ (\pi(\phi)f)(x)
= \int_{\R^\times} \phi(a) |a|^{\frac{s}{2}-1} f(a^{-1}x)\, \frac{da}{|a|}
= \int_{\R^\times} \phi(a) |a|^{\frac{s}{2}-2} f(a^{-1}x)\, da.\]
The antilinear functional
$\alpha$ on $C_c(\R^\times) \subeq \cH_s$,
given by
\begin{equation}
  \label{eq:delfunc}
 \alpha(f) :=  \int_{\R^\times}  \oline{f(x)} |x-1|^{-s}\, dx,
\end{equation}
is formally the scalar product $\la \delta_1, f\ra$ with the
$\delta$-function in $1$. We then have
\begin{align*}
\alpha(\pi(\phi^*)f)
&= \int_{\R^\times}  \int_A \phi(a^{-1}) |a|^{s/2-1} \oline{f(a^{-1}x)}
\, \frac{da}{|a|} |x-1|^{-s}\, dx \\
&= \int_A {\phi(a^{-1})} |a|^{s/2-1}
\int_{\R^\times}  \oline{f(a^{-1}x)} \, |x-1|^{-s}\, dx \frac{da}{|a|}  \\
&= \int_A {\phi(a^{-1})} |a|^{s/2}
\int_{\R^\times}  \oline{f(x)} \, |ax-1|^{-s}\, dx \frac{da}{|a|}  \\
&= \int_A {\phi(a)} |a|^{-s/2}
\int_{\R^\times}  \oline{f(x)} \, |a^{-1}x-1|^{-s}\, dx \frac{da}{|a|}  \\
&= \int_A {\phi(a)} |a|^{s/2}
\int_{\R^\times}  \oline{f(x)} \, |x-a|^{-s}\, dx \frac{da}{|a|}  \\
&= \int_{\R^\times} \int_{\R^\times}
{\phi(a)} |a|^{s/2-1} \oline{f(x)} \, |x-a|^{-s}\, dx\, da.
\end{align*}
Therefore $\alpha \circ \pi(\phi^*)$ can be identified with the
element of $\cH_s$ given by
\[ \big(\alpha \circ \pi(\phi^*)\big)(x) = \phi(x) |x|^{\frac{s}{2}-1}.\]
Since the map
\[ \cD(A) \to \cH_s, \quad \phi \mapsto (x \mapsto \phi(x) |x|^{\frac{s}{2}-1})\]
is continuous, $\alpha$ defines a distribution vector
in $\cH_s^{-\infty}$ with
$\pi^{-\infty}(\phi)\alpha = \alpha \circ \pi(\phi^*)$ for every
$\phi \in \cD(A)$ (cf.~Proposition~\ref{prop:2.4}).  
As $C_c(\R^\times)$ is dense in $\cH_s$, this distribution vector
is cyclic for the representation of $A$ on $\cH_s$.
{}From \eqref{eq:delfunc} we further derive that $\theta \alpha = \alpha$:
\begin{align*}
\la \alpha, \theta  f \ra
&=  \int_\R |x|^{s-2} \oline{f(x^{-1})} |x-1|^{-s}\, dx
=  \int_\R |x|^{s-1} \oline{f(x^{-1})} |x-1|^{-s}\, \frac{dx}{|x|}\\
&=  \int_\R |x|^{1-s} \oline{f(x)} |x^{-1}-1|^{-s}\, \frac{dx}{|x|}
=  \int_\R \oline{f(x)} |1-x|^{-s}\, dx = \la \alpha, f\ra.
\end{align*}
The corresponding distribution on $A$ is given by
\[ \pi^{\alpha}(\phi)
= \la \alpha, \pi^{-\infty}(\phi)\alpha \ra
= \int_\R  \oline{\phi(x)} |x|^{\frac{s}{2}-1} |x-1|^{-s}\, dx
= \int_\R  \oline{\phi(x)} |x|^{\frac{s}{2}} |x-1|^{-s}\, \frac{dx}{|x|},\]
hence represented by the locally integrable function
\[ \gamma(x) := |x|^{\frac{s}{2}} |x-1|^{-s} \quad \mbox{ on } \quad A = \R^\times,\]
which is $\tau$-invariant. On the open subsemigroup
$S = \{ a \in A \: |a| < 1 \}$ we have the expansion
\[ |x|^{s/2}(1-x)^{-s}= \sum_{n=0}^\infty {n-1+s \choose n} |x|^{s/2} x^n.\]
Accordingly, the corresponding representation of $S$ on $\hat\cE$ decomposes as a direct
sum of one-dimensional representations corresponding to the characters
$|x|^{s/2} x^n$, $n \in \N_0$.

Since $\cH_s$ is a space defined by the positive definite
distribution $D = |x|^{-s}$ on $\R$, it also has a natural realization
$\Phi \: \cH_s \to \cD'(\R)$ as a space $\cH_D$ of distributions.
This map is given by
\[ \Phi(f)(\phi) = \int_\R \int_\R \oline{\phi(x)} f(y)|x-y|^{-s}\, dx\, dy\]
(cf.\ Example~\ref{ex:2.3}(a)). In particular, the distribution
$\Phi(f)$ is represented by a continuous function~$\Gamma(f)$. In this sense,
we have
\begin{align*}
 \Phi(\pi(a)f)(x)
&= |a|^{s/2-1} \int_\R f(a^{-1}y)|x-y|^{-s}\, dy
= |a|^{s/2} \int_\R f(y)|x-ay|^{-s}\, dy \\
&= |a|^{-s/2} \int_\R f(y)|\frac{x}{a}-y|^{-s}\, dy
= |a|^{-s/2} \Gamma(f)(x/a),
\end{align*}
so that the corresponding representation of $A$
on $\Phi(C_c(G))$ is given by
\[ (\hat\pi(a)\Gamma(f))(x) := |a|^{-s/2} \Gamma(f)(x/a).\]
\end{ex}

\subsection{Extending distributions}

In this subsection we briefly discuss the existence of reflection
positive extensions of positive definite distributions
on the open subsemigroup $S = (\R_+,+)$ of $\R$ to the whole
real line.

For each real $\alpha > -1$, we obtain a measure
$\mu_\alpha := x^\alpha dx$ on $[0,\infty[$
whose Laplace transform exists and satisfies
$\cL(\mu_\alpha) = c_\alpha x^{-1-\alpha}$.
This shows that the functions $x^{-s}$, $s > 0$, on $S$
are positive definite (cf.\ \eqref{eq:lapform} in Example~\ref{ex:2.3}).

We write $r^{-s} = |x|^{-s}$ for their symmetric extensions
to $\R^\times$. For $s < 1$, $r^{-s}$ is locally integrable,
so that it defines a distribution on $\R$, and by
Lemma~\ref{lem:1.3}, this distribution is positive definite.

For $s \geq 1$, the situation is more complicated.
According to \cite[Thm.~VIII, \S VII.4]{Schw73}, the
measure $|x|^{-s}\, dx$ on $\R^\times$ extends to a distribution
on $\R$. To describe such extensions more explicitly,
one applies the ``finite part'' technique to
the restrictions to $]0,\infty[$ and $]-\infty,0[$
(see \cite[(II.2;26)]{Schw73} for the functions
$r^m$ on $]0,\infty[$). One shows that, for
$\phi \in \cD(\R)$,  the function
$I(\phi,\eps) := \int_{|x|> \eps} \phi(x)|x|^{-s}\, dx$
has an asymptotic expansion for $\eps > 0$ of the form
\[ I(\phi,\eps) = \sum_{k,\ell} a_{k\ell}(\phi) \eps^{-k} (\log \eps)^\ell\]
and, for $s \leq -1$, one defines $\Pf(r^{-s}) \in \cD'(\R)$ by
\[ \Pf(r^{-s})(\phi) := a_{00}(\phi).\]
If $s\not\in -1 - 2\N_0$, then the Fourier transform of
$\Pf(r^{-s})$ is a multiple of $\Pf(r^{s-1})$, and
for $\Pf(r^{-1-2h})$, $h \in \N_0$, it is a sum of a multiples of $r^{2h}$ and
$r^{2h} \log r$ (\cite[(VII.7;13)]{Schw73}).
We conclude that, for each $s \geq  1$, there
exists a polynomial $P$ for which $\cF(\Pf(r^{-s}))(x) + P(4\pi^2 x^2)
\geq 0$. Therefore $\Pf(r^{-s}) + P(-\frac{d^2}{dx^2})\delta_0$
is a reflection positive extension of $r^{-s}$.

\section{A Bochner--Schwartz Theorem for involutive cones} \mlabel{app:a}

In this section we study the convolution
algebra $\cD(S) = C^\infty_c(S,\C)$ of complex-valued
test functions on an open subsemigroup $S$ of a symmetric Lie group
$(G,\tau)$. We assume that $\1 \in \oline S$ and that
$S$ is invariant under the involution $s^\sharp = \tau(s)^{-1}$.
Endowed with the $L^1$-norm, we obtain on
${\cal D}(S)$ the structure of a normed
$*$-algebra $\cA$ with an approximate identity.
We are mainly interested in the case where
$G = V$ is a finite-dimensional real vector space and
$S = \Omega$ is an open convex cone. In this case we show that
the spectrum $\hat\cA$ (=the set of non-zero bounded characters)
of $\cA$ is homeomorphic to a
certain closed convex cone $\hat\Omega \subeq V_\C^*$
parametrizing the bounded $*$-homomorphisms $S \to \C$.
Our main result is a generalization of the Bochner--Schwartz Theorem
(Theorem~\ref{thm:a.7}) characterizing positive definite distributions
$D \in \cD'(S)$ corresponding to contraction representations
of $S$ as Fourier--Laplace transforms of regular
Borel measures on $\hat\Omega$.

\subsection{Some generalities on open subsemigroups}

Let $(G,\tau)$ be a symmetric Lie group and
$S \subeq G$ be an open subsemigroup with
$\1 \in \oline S$ which is invariant under the involution
$s^\sharp = \tau(s)^{-1}$.
Then ${\cal D}(S)$ is a convolution subalgebra
of ${\cal D}(G)$ with an approximate identity and involution
$\phi^\sharp := \phi^* \circ \tau$.
We define a representation of ${\cal D}(S)$ on the space
${\cal D}'(S)$ by
\[ (\phi * D)(\psi) := D(\phi^\sharp * \psi). \]

The following proposition generalizes the classical
Dixmier--Malliavin Theorem (\cite[Thm.~3.1]{DM78}) to open
subsemigroups of Lie groups.

\begin{prop} \mlabel{prop:2.5} Let $S \subeq G$ be an open subsemigroup with
$\1 \in \oline S$. Then every test function
$\phi \in \cD(S)$ is a sum of products
$\alpha * \beta$ with $\alpha, \beta \in \cD(S)$.
\end{prop}

\begin{prf} Let $\phi \in \cD(S)$.
Then $\supp(\phi)$ is a compact subset of $S$, so that there exists a
$\1$-neighborhood $U \subeq G$ with $U \supp(\phi) \subeq S$.
Pick $u \in U \cap S^{-1}$ and let $V \subeq G$ be a $\1$-neighborhood in
$G$ with $u^{-1}V \subeq S$. Then
$\phi \circ \lambda_u^{-1} = \delta_u*\phi \in \cD(S)$.

According to \cite[Thm.~3.1]{DM78}, $\delta_u * \phi$ is a finite
sum of functions of the form
$\alpha * \beta$ with
$\supp(\alpha) \subeq V$ and $\supp(\beta) \subeq \supp(\delta_u * \phi)
= u \supp(\phi)$. Then $\phi$ is a finite sum of functions of the
form $\delta_{u^{-1}} * \alpha * \beta$, with
\[ \supp(\delta_{u^{-1}} * \alpha) =
u^{-1} \supp(\alpha) \subeq u^{-1}V \subeq S \]
and
$\supp(\beta) \subeq u \supp(\phi) \subeq S.$
\end{prf}

\begin{lem}
  \mlabel{lem:3.1}
For each $D \in {\cal D}'(S)$ and $\phi \in {\cal D}(S)$ the distribution
$\phi * D$ is represented by the smooth function
\[ (\phi * D)(x) := D((\phi \circ \lambda_{\tau(x)})^\sharp)\]
on $S$ which extends smoothly to an open neighborhood of $\oline S$.
\end{lem}

\begin{prf} The first assertion follows easily
from \cite[Prop.~A.2.4.1]{Wa72}
if we take into account that we defined distributions as antilinear
functionals on $\cD(G)$.

For the second assertion we note that the function
\[ \Phi\: G \to \cD(G), \quad x \mapsto \phi \circ \lambda_{\tau(x)} \]
is smooth. As $\supp(\Phi(x)) = x^\sharp \supp(\phi)$ and
$\supp(\phi)$ is a compact subset of $S$,
the set
$S_\phi := \{ x \in G \: x^\sharp \supp(\phi) \subeq S \}$ is an open
neighborhood of $\oline S$, and the function
\[  \Phi \:
\{ x \in G \: x^\sharp \supp(\phi) \subeq S \} \to {\cal D}(S),\]
is smooth. This  completes the proof.
\end{prf}

\begin{lem}
  \mlabel{lem:3.2}
If $\chi \in {\cal D}'(S)$ is a non-zero homomorphism of
$*$-algebras, then
$\chi$ is represented by a smooth homomorphism
$S \to (\C, \cdot)$ of involutive semigroups.
If, in addition, for every $\1$-neighborhood $U$ in $G$ the set
$S \cap U$ generates $S$, then $\chi$ has no zeros.
\end{lem}

\begin{prf}
Since $\chi$ is non-zero, there exists $\phi \in {\cal D}(S)$
with $\chi(\phi)\not=0$.
Then
\[ (\phi * \chi)(\psi) = \chi(\phi^\sharp * \psi) = \chi(\psi)
\oline{\chi(\phi)} \]
implies that
\begin{equation}
  \label{eq:a.2}
\chi = \oline{\chi(\phi)}^{-1} (\phi * \chi),
\end{equation}
and hence that $\chi$ is represented by a smooth function 
which we also denote by $\chi$ (Lemma~\ref{lem:3.1}).

For $\psi \in \cD(S)$ the relation
\begin{align*}
\int_S \phi(s)\oline{\chi(s)}\, d\mu_G(s)
&= \oline{\chi(\phi)}
=\chi(\phi^\sharp)\\
&=\int_S \phi(s^\sharp)\Delta(s^{-1})\chi(s)\, d\mu_G(s)
=\int_S \phi(s)\chi(s^\sharp)\, d\mu_G(s)
\end{align*}
implies that $\chi(s^\sharp) = \oline{\chi(s)}$.

Using \eqref{eq:a.2}, we obtain the relation
\begin{equation}
  \label{eq:mulrel}
\oline{\chi(\phi)}\chi(x)
= \chi((\phi \circ \lambda_{\tau(x)})^\sharp)
= \oline{\chi(\phi \circ \lambda_{\tau(x)})}.
\end{equation}
Now we let $\phi_n$ be a sequence of test functions with total integral $1$
such that $\supp(\phi_n)$ converges to the point $y \in S$.
Passing to the limit in \eqref{eq:mulrel} now leads to
\[ \chi(y^\sharp)\chi(x)   = \oline{\chi(y)}\chi(x) =  \oline{\chi(x^\sharp y)}
= \chi(y^\sharp x).\]
Therefore $\chi \: S \to (\C, \cdot)$ is a homomorphism.

It remains to show that $\chi(x)\not= 0$ for all $x \in S$
if $S$ is generated by
every neighborhood of $\1$, intersected with $S$. Pick $x \in S$ with
$\chi(x) \not= 0$. Then the open set $xS^{-1}$ is a neighborhood of $\1$
and for
$z \in S \cap xS^{-1}$ the relation
$0\not= \chi(x) = \chi(z) \chi(z^{-1}x)$ implies that
$\chi(z)\not= 0$. Since $S$ is generated by $S \cap xS^{-1}$, the assertion follows.
\end{prf}

\begin{lem}
  \mlabel{lem:3.3} The norm
$$ \|\phi\|_1 := \int_S |\phi(x)|\, d\mu_G(x) $$
on $\cD(S)$
is submultiplicative, so that $({\cal D}(S), \|\cdot\|_1)$ is a normed algebra
and $\phi^\sharp := \phi^* \circ \tau$ defines an isometric involution
on ${\cal D}(S)$.
\end{lem}

\begin{prf} We know already from Definition~\ref{def:2.1} that
$(\cD(S), *, \sharp)$ is an involutive algebra.
The submultiplicativity of $\|\cdot\|_1$
and $\|\phi^\sharp\|_ = \|\phi\|_1$ follow immediately from the
corresponding properties of $L^1(G)$ because the involution
$\tau$ preserves the Haar measure on~$G$.
\end{prf}

\begin{defn} \mlabel{def:a.3}
We write $\hat S$ for the set of all continuous
homomorphisms $\chi \: S \to \C$ satisfying
\[ \chi(s^\sharp) = \oline{\chi(s)}\quad \mbox{ and } \quad
|\chi(s)| \leq 1 \quad \mbox{ for }  \quad s\in S.\]

It is easy to see that any $\chi \in \hat S$ defines a
$*$-homomorphism
\[ \cD(S) \to \C,\quad \phi \mapsto \int_S \phi(s)\chi(s)\, d\mu_G(s),\]
hence is smooth by Lemma~\ref{lem:3.2}.
\end{defn}

\begin{rem}
  \mlabel{rem:3.4}
Let $D \in {\cal D}'(S)$ be a positive definite distribution and
$(\pi_D, {\cal H}_D^0)$ be the corresponding representation, where
${\cal H}_D^0 = {\cal D}(S) * D \subeq \cD'(S)$
and $\pi_D(\phi)T = \phi * T$.

(a) We assume, in addition, that $D$ is {\it $1$-bounded}
in the sense that
$$ \|D(\phi* \psi * \phi^*)| \leq \|\psi\|_1 D(\phi*\phi^*)
\quad \mbox{ for }  \quad \phi, \psi \in {\cal D}(S). $$
Then the representation of ${\cal D}(S)$ on ${\cal H}_D^0$
extends to a representation on
the reproducing kernel Hilbert space
${\cal H}_D \subeq {\cal D}'(S)$ (cf.\ \cite[Th.\ III.1.19]{Ne00}
and Section~\ref{sec:2}).
If $D' \in {\cal D}'(S)$ satisfies
${\cal D}(S) * D' = \{0\}$, then the existence of an approximate identity in $S$ implies
that $D' = 0$. Therefore the representation of ${\cal D}(S)$ on ${\cal H}_D$
is non-degenerate.

(b) Suppose, in addition, that $S$ is commutative.
We claim that $\pi_D(\phi) = 0$ implies that
$\phi * D = 0$. In fact, $\pi_D(\phi) = 0$ leads
to
\[ 0 = \pi_D(\phi)(\psi * D)= \phi * \psi * D = \psi * \phi * D
\quad \mbox{ for }\quad \psi \in \cD(\Omega), \]
so that
\[ D(\phi^\sharp * \psi^\sharp * \xi) = 0\quad \mbox{ for } \quad
\xi, \psi \in \cD(\Omega),\]
i.e., $\phi *D$ vanishes on
$\cD(S) * \cD(S) = \cD(S)$
(Proposition~\ref{prop:2.5}). This means that $\phi * D = 0$.
\end{rem}

\subsection{Open convex cones}

We now turn to the special case where $G = V$ is a
finite-dimensional real vector space and
$\tau \in \GL(V)$ is an involution.
Let $\Omega \subeq V$ be an open convex cone
invariant under the involution
$s^\sharp := - \tau(s)$. Then $(\Omega, \sharp)$ is an involutive
semigroup and the convolution algebra $\cD(\Omega)$ of test functions
on $\Omega$ is an involutive algebra with respect to the
involution $\phi^\sharp(s) := \oline{\phi(s^\sharp)}$.

Since $\Omega$ is generated by any intersection
$U \cap \Omega$, where $U \subeq V$ is a $0$-neighborhood,
any non-zero $\chi \in \hat\Omega$ (the set of bounded
$*$-homomorphisms into $(\C,\cdot)$) maps into $\C^\times$.
Hence $\chi$ extends to a $*$-homomorphism
$(V,\sharp) \to \C^\times$, which means that
$\chi = e^{-\alpha}$ for a linear functional
$\alpha  \: V \to \C$ satisfying $\alpha^\sharp := - \oline\alpha \circ \tau
= \alpha$. The set of all these functionals
can be identified with the dual space of
\[ V_c := \{ x + i y \in V_\C  \:  x^\sharp = x, y^\sharp = - y\}.\]
In the following we therefore identify $\hat\Omega$
with the closed convex cone
\[ \hat\Omega = \{ \alpha \in V_c^* \: (\forall
x = x^\sharp \in \Omega)\ \alpha(x) \geq 0\}\]
in the real vector space $V_c^*$.

\begin{exs} \mlabel{ex:a.1}
(a) $\tau = -\id_V$, $\sharp = \id_V$, and $\Omega$ any open convex cone.
Then $V_c = V$ and $\hat\Omega = \Omega^\star$ is the dual cone.

(b) $\tau = \id_V$ and $V = \Omega$.
Then $V_c = iV$ and $\hat\Omega = V_c^* = iV^*$.

(c) $V = \R^{n}$, $\tau(x_0, x_1,\ldots, x_{n-1}) =
(-x_0, x_1, \ldots, x_{n-1})$ and
\[ \Omega  := \R^n_+ = \{ x \in \R^{n} \: x_0 > 0\} \]
an open half space. Then
\[ V_c = \R e_0 + i \sum_{j > 0} \R e_j
\quad \mbox{ and } \quad
\hat\Omega = \{ \alpha \in V_c^* \: \alpha_0 \geq 0\}.\]
\end{exs}

\begin{rem} \mlabel{rem:4.8}
For the case $\Omega = V$ with $\tau(v) = v$, the
Bochner--Schwartz Theorem (\cite[Thm.~XVIII, \S VII.9]{Schw73})
asserts that $D \in \cD'(V)$ is positive definite if and only if
it is tempered and its Fourier transform $\hat D$ is a positive measure.
We then have an embedding $\cH_D \subeq \cS'(V)$ (cf.\ Example~\ref{ex:2.7}).
\end{rem}

\begin{defn} (a) Let $\cA := C^*(\cD(\Omega),\|\cdot\|_1)$ be
the enveloping $C^*$-algebra of the normed involutive algebra
$(\cD(\Omega), \|\cdot\|_1)$ and $\eta \: \cD(\Omega) \to \cA$ be the
canonical map.

(b) For each
$\phi \in \cD(\Omega)$, the  {\it Fourier--Laplace transform}
\[ \hat\phi \: \hat\Omega \to \C, \quad
\hat\phi(\alpha) := \int_\Omega \phi(x)e^{-\alpha(x)}\, dx \]
is a continuous function on $\hat\Omega$. It corresponds to evaluation
of the character $e^{-\alpha}$ of $\hat\Omega$ on $\phi$. In
particular, we have
\[ (\phi * \psi)\,\hat{} = \hat\phi \cdot \hat\psi
\quad \mbox{ and } \quad
\hat{\phi^\sharp} = \oline{\hat\phi} \quad \mbox{ for }
\quad \phi, \psi \in \cD(\Omega).\]
\end{defn}

\begin{prop} \mlabel{prop:spec} The Fourier--Laplace transform
defines a homeomorphism
\[ \Gamma \: \hat\Omega \to \hat\cA, \quad
\Gamma(\alpha)(\eta(\phi)) := \hat\phi(\alpha) =
\int_\Omega \phi(s)e^{-\alpha(s)}\, ds, \quad \phi \in \cD(\Omega),\]
where $\hat\Omega$ carries the topology induced from $V_c^*$ and
$\hat\cA \subeq \cA'$ the weak-$*$-topology.
\end{prop}

\begin{prf}  By definition, the spectrum $\hat\cA$ of $\cA$ is the set
of all non-zero one-dimensional representations of
$\cA$, which is the same as the set of all non-zero $*$-homomorphisms
$\chi \: \cD(\Omega) \to \C$ with $\|\chi\|_1 \leq 1$.
In view of Lemma~\ref{lem:3.2}, any such character is represented
by a $*$-homomorphism $\beta \: \Omega \to \C^\times$ and the
condition $\|\chi\|_1 \leq 1$ corresponds to the boundedness of
$\beta$. We conclude that $\Gamma$ is bijective.

For $\phi \in \cD(\Omega)$ its Fourier--Laplace transform
$\hat\phi \: \hat\Omega \to \C$ is continuous and
since $\eta(\cD(\Omega))$ is dense in $\cA$ and $\hat\cA$ is bounded,
it follows that $\alpha \mapsto \Gamma(\alpha)(A)$
is continuous for every $A \in \cA$, i.e., that
$\Gamma$ is continuous. It remains to show that its inverse is
also continuous.

To this end, we first observe that the translation action of $\Omega$ on
$\cD(\Omega)$
\[ (s.\phi)(x) := \phi(x-s) \]
defines multipliers of the $C^*$-algebra $\cA$,
which leads to a homomorphism
\[ \eta \: \Omega \to M(\cA) \cong C_b(\hat\cA) \]
(\cite[Ex.~12.1.1(b)]{Bl98}). More concretely,
\[ (s.\phi)\,\,\hat{}(\alpha)
= \int_\Omega \phi(x-s)e^{-\alpha(x)}\, dx
= \int_{\Omega-s} \phi(x)e^{-\alpha(x+s)}\, dx
= e^{-\alpha(s)}\hat\phi(\alpha).\]
For every $A \in \cA$, the
map $\Omega \to \cA, x \mapsto \eta(x)A$ is continuous, i.e.,
$\eta$ is continuous with respect to the so-called strict
topology on $M(\cA)$. Since $\eta(\Omega)$ is bounded
and this topology coincides on bounded subsets
with the compact open topology (\cite[Lemma A.1]{GrN09}),
the map $\eta \: \Omega \to C_b(\hat\cA)$ is continuous with
respect to the compact open topology on $C_b(\hat\cA)$.
This in turn means that the map
\[ \Omega \times \hat\cA \to \C, \quad (s,\chi) \mapsto
\chi(\eta(s)) \]
is continuous (\cite[Thm.~VII.2.4]{Br93}), and hence that the map
\[ \hat\eta \: \hat\cA \to \hat\Omega,
\quad \hat\eta(\chi)(s) := \chi(\eta(s)) \]
is continuous if $\hat\Omega$ carries the topology
of uniform convergence on compact subsets of~$\Omega$. It is easy to see that this topology coincides
with the subspace topology inherited from $V_c^*$.
This proves that $\Gamma$ is a homeomorphism.
\end{prf}

\begin{thm} \mlabel{thm:a.7} {\rm(Generalized Bochner--Schwartz Theorem)}
Let $\tau$ be an involution on the finite-dimensional vector space
$V$ and $\Omega \subeq V$ be an open convex cone invariant under the
involution $v \mapsto v^\sharp := -\tau(v)$.
If $D \in \cD'(\Omega)$ is a positive definite $1$-bounded
distribution, then the following assertions hold: 
\begin{description}
\item[\rm(a)] There exists a unique positive Radon measure
$\mu$ on the closed convex cone $\hat\Omega$ with
\[ D(\phi) = \int_{\hat\Omega} \oline{\hat\phi(\alpha)}\, d\mu(\alpha)
\quad \mbox{ for } \quad \phi \in \cD(\Omega),\]
resp.,
\[ D  = \int_{\hat\Omega} e_{-\oline\alpha}\,  d\mu(\alpha) \quad \mbox{ for }  \quad 
e_{-\alpha}(x) = e^{-\alpha(x)}, \quad x \in \Omega.\] 
\item[\rm(b)] There exists a unitary map 
$\Gamma \: \cH_D \to L^2(\hat\Omega, \mu)$ mapping 
$\phi * D$ to $\hat\phi$ and intertwining the contraction representation 
of $\Omega$ on 
$\cH_D \subeq \cD'(\Omega)$ with the multiplication representation 
$\pi_\mu(s)f = e_{-s} f$. 
\end{description}
\end{thm}

\begin{prf}  (a) The representation
$(\pi_D, \cH_D^0)$ of $\cD(\Omega)$ leads to a representation
$\tilde\pi_D \: \cA \to B(\cH_D)$ with $\tilde\pi_D  \circ \eta
= \pi_D$ and $\pi_D(\phi) = 0$ implies that
$\phi * D = 0$ (Remark~\ref{rem:3.4}).
We conclude that the map
\[ \gamma \: \cD(\Omega) \to \cH_D, \quad \phi \mapsto \phi * D \]
factors through a linear map
\[ \oline\gamma \: \pi_D(\cD(\Omega)) \to \cH_D, \quad
\pi_D(\phi) \mapsto \phi * D. \]
Hence
\[ K(\pi_D(\phi),\pi_D(\psi))
:= \la \phi * D, \psi * D \ra_{\cH_D}
= D(\phi^\sharp * \psi) \]
is a positive definite sesquilinear kernel on
the dense subalgebra $\cB := \pi_D(\cD(\Omega))$
of the $C^*$-algebra $\tilde\pi_D(\cA)$.

{}From Proposition~\ref{prop:2.5} we know that
$\cD(\Omega) * \cD(\Omega) = \cD(\Omega)$, which also implies
that $\cD(\Omega) * \cD(\Omega) * \cD(\Omega) = \cD(\Omega)$. Therefore
\cite[Prop.~II.4.13]{Ne00} shows that the kernel
$K$ is non-degenerate because all totally degenerate kernels on the
$*$-algebra $\cD(\Omega)$ vanish. By assumption, the kernel
$K$ is $\|\cdot\|_1$-bounded. Now
the abstract Plancherel Theorem \cite[Thm.~VI.1.6]{Ne00} 
implies the existence of a unique positive Borel measure on
$\hat \cA \cong \hat\Omega$ (Proposition~\ref{prop:spec})
with
\[ D(\phi^\sharp * \psi)
= \int_{\hat\Omega} \hat\phi(\alpha)\oline{\hat\psi(\alpha)}\, d\mu(\alpha)
= \int_{\hat\Omega} (\phi * \psi^\sharp)\,\hat{}(\alpha)
\, d\mu(\alpha).\]
Applying Proposition~\ref{prop:2.5} again, it follows that
\[ D(\phi) = \int_{\hat\Omega} \oline{\hat\phi(\alpha)}\, d\mu(\alpha)
\quad \mbox{ for } \quad \phi \in \cD(S).\]

(b) follows from \cite[Thm.~VI.1.6]{Ne00}. 
\end{prf}

For a distribution $D_f(\phi) = \int_\Omega \oline\phi(x) f(x)\, dx$
given by a positive definite function $f$ on $\Omega$,
we apply the integral representation from the
preceding theorem to $\delta$-sequences in the point
$s \in \Omega$ to obtain
\[ f(s)
= \int_{\hat\Omega} \oline{e^{-\alpha(s)}}\, d\mu(\alpha)
= \int_{\hat\Omega} e^{-\oline\alpha(s)}\, d\mu(\alpha).\]
This is the integral representation obtained by
Shucker in \cite[Thm.~5]{Sh84} for continuous positive definite
functions $r \: \Omega \to \C$
on convex domains $\Omega \subeq V$ satisfying
$\Omega + V^\tau = \Omega$ for which the kernel
\[ K(z,w) = r\Big(\frac{z + w^\sharp}{2}\Big) \]
is positive definite. For the case of convex cones, the preceding
theorem extends Shucker's Theorem to distributions.

\begin{ex} (cf.\ Example~\ref{ex:a.1}) (a) 
For an open convex cone $\Omega \subeq V$ 
with the involution $\tau= -\id_V$, we have $\hat\Omega = \Omega^\star$ 
and we obtain an integral representation of the form
\[ D  = \int_{\Omega^\star} e_{-\alpha}\,  d\mu(\alpha), \quad 
e_{-\alpha}(x) = e^{-\alpha(x)}, \quad x \in \Omega.\] 

(b) For $\tau = \id_V$ and $V = \Omega$ we have 
$\hat\Omega = iV^*$ and we recover the 
classical Bochner--Schwartz Theorem asserting that every 
positive definite distribution on the group $(V,+)$ 
is the Fourier transform 
\[ D  = \int_{V^*} e_{i\alpha}\,  d\mu(\alpha) \] 
of a measure on $V^*$ (cf.\ Remark~\ref{rem:4.8}). 
 
(c) For $V = \R^{n}$, $\tau(x_0, x_1,\ldots, x_{n-1}) =
(-x_0, x_1, \ldots, x_{n-1})$ and the open half space 
$\Omega  = \R^n_+$ 
we obtain an integral representation 
\[ D(x_0,x')  = \int_0^\infty \int_{\partial\Omega} e^{-\lambda x_0 + i \la y, x'\ra}
\,  d\mu(\lambda,y).\] 
\end{ex}

In the case where $\tau = - \id_V$, we have the following 
sharper version of the Bochner--Schwartz Theorem. It implies in particular that 
positive definite distributions are functions. 

\begin{thm}
  \mlabel{thm:4.1} Suppose that $\tau = - \id_V$ and that 
$\Omega \subeq V$ is an open convex cone not containing affine lines. 
Then any $1$-bounded positive definite distribution 
$D \in \cD'(\Omega)$ is represented by an analytic function on $\Omega$, 
also denoted $D$,  
and there exists a unique Radon measure $\mu$ on the dual cone 
$\Omega^\star$ with 
\[  D(x) = {\cal L}(\mu)(x) = \int_{\Omega^\star} e^{-\alpha(x)}\, d\mu(\alpha). \]
\end{thm}

\begin{prf} First we use Theorem~\ref{thm:a.7} to obtain an 
integral representation 
$$ D = \int_{\Omega^\star} e_{-\alpha}\, d\mu(\alpha), $$ 
where $\mu$ is a positive Radon measure on the dual cone $\Omega^\star \subeq V^*$. 
For $\phi \in {\cal D}(\Omega)$ this leads to 
\[  D(\phi) = \int_{\Omega^\star} \la e_\alpha, \phi \ra\, d\mu(\alpha), 
\quad \mbox{ where } \quad 
\la e_\alpha, \phi \ra = \int_\Omega e^{-\alpha(x)} \phi(x)\, dx =:{\cal L}(\phi)(\alpha) \]
and $dx$ stands for a Haar measure on $V$. 
It follows that, for every $\phi \in {\cal D}(\Omega)$, the Laplace transform 
${\cal L}(\phi)$, viewed as a function on $\Omega^\star$, is integrable with respect to $\mu$.

Let $x \in \Omega$. Since $x - \Omega$ is a $0$-neighborhood and $0 \in \oline \Omega$, 
there exists a non-negative $\phi \in {\cal D}(\Omega)$ with $\int_\Omega \phi(x)dx = 1$ and 
$\supp(\phi) \subeq x - \Omega$. Then we have $e^{-\alpha(x)} \leq e^{-\alpha(y)}$ for all 
$y \in \supp(\phi)$ and $\alpha \in \Omega^\star$, so that we also get 
$e^{-\alpha(x)} \leq {\cal L}(\phi)(\alpha).$
Therefore the integrability of ${\cal L}(\phi)$ implies the integrability of 
the function $e_{-x}(\alpha) := e^{-\alpha(x)}$ for all $x \in \Omega$. Now 
\[ {\cal L}(\mu)(x) := \int_{\Omega^\star} e^{-\alpha(x)}\, d\mu(\alpha) \] 
exists as a function on $\Omega$. According to \cite[Prop.~V.4.6]{Ne00}, 
the function ${\cal L}(\mu)$ has a 
holomorphic extension to the tube domain $\Omega + i V \subeq V_\C$, hence is in particular analytic. 
That $D$ is represented by the function ${\cal L}(\mu)$ follows by Fubini's Theorem 
which applies because all functions 
$(x,\alpha) \mapsto e^{-\alpha(x)}\phi(x)$ are integrable with respect to 
the product measure $dx \otimes \mu$ on $\Omega \times \Omega^\star$.
\end{prf}

\section{The complementary series of the conformal group} \mlabel{sec:4}

In this section we discuss the complementary series $(\pi_s, \cH_s)_{0 < s < n}$
of irreducible unitary representations of the group
$\OO_{1,n+1}(\R)$ which acts by partially defined conformal transformations
on $\R^n$ and by everywhere defined maps on its conformal completion $\bS^n$ 
(cf.\ \cite{vD09}). These representations are reflection positive for 
several involutions, resp., semigroups, and this leads to a remarkable 
connection of the distributions $\|x\|^{-s}$ on $\R^n$ with positive definite 
kernels on the open unit ball. 

\subsection{The conformal group on $\R^n$}

We describe the {\it conformal completion of $\R^n$} by
the stereographic map
\begin{equation}
  \label{eq:stereo}
\eta \: \R^n \to \bS^n, \quad
\eta(x) = \Big( \frac{1-\|x\|^2}{1+\|x\|^2},
\frac{2 x}{1 + \|x\|^2}\Big)
\end{equation}
whose image is the complement of the point
$-e_0 \in \bS^n$.

To obtain an action of the conformal group on the sphere,
we embed $\bS^n$ into the set of isotropic vectors of the
$(n+2)$-dimensional Minkowski space
\[ \zeta \: \bS^n \to \R^{n+2}, \quad
\zeta(y) = (1,y).\]
The image of this map intersects each isotropic line in $\R^{n+2}$
exactly once, so that it leads to a diffeomorphism of $\bS^n$ with the
projective quadric of isotropic lines in $\bP(\R^{n+2})$.

The group $\OO_{1,n+1}(\R)$ acts on $\R^{n+2}$ preserving the Lorentzian
form. Writing elements of $\R^{n+2}$ as pairs
$z= (z_0, z')$ with $z_0 \in \R, z' \in \R^{n+1}$, we find for
$g = \pmat{ a & b \\ c & d} \in \OO_{1,n+1}(\R)$ the relation
\begin{equation}
  \label{eq:confact-}
gz = (az_0 + \la b,z'\ra, cz_0 + dz').
\end{equation}
If $z\not=0$ satisfies $[z,z]=0$, then $z_0 \not=0$. This leads to
$a z_0 + \la b,z'\ra\not=0$. Therefore the induced action of
$\OO_{1,n+1}(\R)$ on
\[ \bS^n \cong \{ z = (1,x) \: \|x\| =1 \}
\cong \{ 0 \not=z \in \R^{n+2} \: [z,z]= 0\}/\R^\times \]
can be written as
\begin{equation}
  \label{eq:confact}
g.x = (a + \la b,x\ra)^{-1}(c + dx), \quad x \in \bS^n \subeq \R^{n+1}.
\end{equation}
In terms of the realization of $\R^n$ as the complement
$\eta(\R^n)$ of the point
$-e_0$ in $\bS^n$, this action of $\OO_{1,n+1}(\R)$ on $\bS^n$ yields an
``action'' on $\R^n$ by rational maps. For $n > 2$, these maps
exhaust all conformal maps on open subsets of $\R^n$
(\cite[Thm.~1.9]{Scho97}).

The stabilizer of the point $-e_0 \in \bS^n$ ``at infinity'',
which is represented by the element $(1,-e_0) \in \R^{n+2}$
 is the stabilizer $P \subeq \OO_{1,n+1}(\R)$ of the isotropic line \break
$\R(1,-e_0) \subeq \R^{n+2}$. This group is a
parabolic subgroup of
$\OO_{1,n+1}(\R)$, which acts by affine conformal maps on $\R^n$.
The latter description easily shows that $P$ is isomorphic
to
\begin{equation}
  \label{eq:affconf}
\Aff_c(\R^n) := \R^n \rtimes (\R^\times \times \OO_n(\R)),
\end{equation}
the group of affine conformal isomorphism of $\R^n$.
In the standard notation for parabolic subgroups, we have
$P = NAM$ with $N \cong \R^n$, acting by translations,
$A \cong \R^\times$ acting by homotheties, and
$M \cong \OO_n(\R)$, the centralizer of $A$ in the maximal
compact subgroup $K \cong \OO_{n+1}(\R)$ of $\OO_{1,n+1}(\R)$.

\begin{rem} \mlabel{rem:4.1}
The group $\OO_{1,n+1}(\R)$ does not act faithfully
on $\bS^n$ because the element $-\1$ acts trivially.
However, the index $2$ subgroup
\[ \OO_{1,n+1}^+(\R) := \{ g \in \OO_{1,n+1}(\R) \: g_{00} > 0\} \]
of those Lorentz transformations preserving the future light cone
\begin{equation}
  \label{eq:lightcon}
\Omega := \{ z \in \R^{n+2} \: [z,z] > 0, z_0 > 0\}
\end{equation}
(cf.\ \eqref{eq:lormet} above),
acts faithfully. Note that $g_{00} = [ g e_0, e_0] \not=0$ follows
from $[g e_0, g e_0] = [e_0, e_0] = 1$.
\end{rem}

\subsection{Canonical kernels}

In this section we discuss the canonical kernel 
$1- \la x,y\ra$ on
$\bS^n$. This kernel is projectively invariant under
$\OO_{1,n+1}(\R)$, and the corresponding cocycle
$J \: \OO_{1,n+1}(\R) \to C^\infty(\bS^n,\R^\times)$
is easily described in terms of the conformal structure.

On the light cone $\Omega \subeq \R^{n+2}$ (cf.\ \eqref{eq:lightcon}),
we have $[x,y]> 0$ for $x,y \in \Omega$, which leads to
the canonical kernel $[x,y]^{-1}$ invariant under the
action of $\OO_{1,n+1}(\R)$.  Restricting to the sphere $\bS^n$ in the
boundary of $\Omega$, we find
\[ [(1,x), (1,y)] = 1 - \la x, y \ra, \]
which leads to the singular kernel $(1 - \la x, y \ra)^{-1}$ on $\bS^n$.

\begin{rem} \mlabel{rem:4.2}
(a) For $g \in \OO_{1,n+1}(\R)$ and $x,y \in \bS^n$, we have
  \begin{align*}
& 1 - \la x,y\ra
= [(1,x),(1,y)] = [g(1,x), g(1,y)] \\
&= [(a + \la b,x\ra, c + dx), (a + \la b,y\ra, c + dy)] 
= (a + \la b,x\ra)(a + \la b,y\ra)[(1, g.x), (1,g.y)]\\
&= (a + \la b,x\ra)(a + \la b,y\ra) (1 - \la g.x, g.y\ra).
  \end{align*}
In view of 
\[ \|x-y\|^2 = 2(1 - \la x,y\ra) \quad \mbox{ for }\quad 
x,y \in \bS^n,\] 
this in turn leads to
\begin{equation}
  \label{eq:squeeze}
\|g.x - g.y\|^2
= (a + \la b,x\ra)^{-1}(a + \la b,y\ra)^{-1}\| x- y\|^2,
\end{equation}
and hence to
\[ \| \dd g(x)v\| = |a + \la b,x\ra|^{-1} \|v\| \quad \mbox{ for }
\quad v \in T_x(\bS^n).\]
The function 
\[ J_g(x) = |a + \la b,x\ra|^{-1} = \|\dd  g(x)\| \]
is called the {\it conformal factor} of $g$. For $g \in \OO_{1,n}^+(\R)$, the
number $a + \la b,x\ra$ is always positive because
it is the $0$-component of $g(1,x)$ (cf.\ Remark~\ref{rem:4.1}).
The Chain Rule implies the cocycle relation
\begin{equation}
  \label{eq:cocyc}
J_{g_1 g_2}(x) = J_{g_1}(g_2.x) J_{g_2}(x),
\end{equation}
and from above we derive for 
$Q_{-1}(x,y) := 1-\la x, y\ra$ the transformation formula 
\begin{equation}
  \label{eq:kertrafo} 
Q_{-1}(g.x, g.y) = J_g(x) Q_{-1}(x,y) J_g(y) \quad \mbox{ for } \quad 
x,y \in \bS^n, g \in \OO_{1,n}^+(\R).
\end{equation}

(b) The stereographic map $\eta \: \R^n \to \bS^n$ from \eqref{eq:stereo}
satisfies
\[ 1 - \la \eta(x), \eta(y)\ra
= 1- \frac{(1-\|x\|^2)(1- \|y\|^2)+ 4 \la x,y \ra}
{(1 + \|x\|^2)(1 + \|y\|^2)}
= 2\frac{\|x-y\|^2}{(1 + \|x\|^2)(1 + \|y\|^2)}, \]
so that  the same argument as above, using
$\|\eta(x)-\eta(y)\|^2 = 2 - 2 \la \eta(x), \eta(y)\ra$,
implies that its conformal
factor, as a map $\R^n \to \bS^n$, is given by
\begin{equation}
  \label{eq:eta-fact}
J_\eta(x) = \frac{2}{1 + \|x\|^2}.
\end{equation}
This implies in particular that, up to a normalizing constant,
\begin{equation}
  \label{eq:invmeas}
(\eta^{-1})_* \dd\mu_\bS
= \frac{2^n}{(1 + \|x\|^2)^n} \, dx
\end{equation}
is the pullback of the surface measure $\mu_\bS$ on $\bS^n$ to $\R^n$.

(c) In view of (b), the pullback of the kernel
$Q(x,y) = (1 - \la x,y\ra)^{-s/2}$ on $\bS^n$ by $\eta$ is 
the kernel
 \begin{equation}
   \label{eq:euc-kern}
K(x,y) := 2^{-s/2}(1 + \|x\|^2)^{s/2} \|x-y\|^{-s} (1 + \|y\|^2)^{s/2}
 \end{equation}
on $\R^n$.

(d) From \eqref{eq:squeeze} it follows in particular, that the cross ratio
\[ \frac{\|y-a\|}{\|y-b\|}\frac{\|x-b\|}{\|x-a\|} \]
is invariant under the action of $\OO_{1,n+1}(\R)$ on the set
of pairwise distinct $4$-tuples in~$\bS^n$.
\end{rem}

\begin{rem} \mlabel{rem:4.3} (A conformal Cayley transform)
We consider the involution $c \in \OO_{1,n+1}^+(\R)$ defined by
\[ ce_0 = e_0,\  ce_1 = e_2,\ c e_2 = e_1
\quad \mbox{ and } \quad ce_j = e_j \quad \mbox{ for } \quad 2 < j \leq n+1.\]
With respect to the action defined in \eqref{eq:confact}, we then have 
\[ c\eta(x) = \frac{1}{1 +\|x\|^2}(2x_0, 1 - \|x\|^2, 2x_1, \ldots,
2 x_{n-1}) = \eta(\phi(x))\]
for the map
\[ \phi(x) = \Big(
\frac{1 - \|x\|^2}{1 + \|x\|^2 + 2 x_0},
\frac{2x_1}{1 + \|x\|^2 + 2 x_0}, \ldots,
\frac{2x_{n-1}}{1 + \|x\|^2 + 2 x_0}\Big).\]
This formula implies that $x_0 > 0$ is equivalent to
$\|\phi(x)\| < 1$ and
$\phi(-e_0) = \infty$. Therefore $\phi$ is an involutive conformal map on
$\R^n$ mapping the open half space $\R^n_+$ onto the open
unit ball $\cD \subeq \R^n$.
For $x_0 = 0$ we have $\phi(x) \in \bS^{n-1}$ and
\[ \phi(0,x) = \Big(
\frac{1 - \|x\|^2}{1 + \|x\|^2},
\frac{2x_1}{1 + \|x\|^2}, \ldots,
\frac{2x_{n-1}}{1 + \|x\|^2}\Big)\]
is the stereographic map $\R^{n-1} \into \bS^{n-1}$ whose image is
the complement of $-e_0$.
\end{rem}

\subsection{A Hilbert space of measures}

We now assume that $0 < s < n$, so that
$\|x\|^{-s}$ is locally integrable on $\R^n$.
{}From the positive definiteness of the distribution $\|x-y\|^{-s}$
(Lemma~\ref{lem:1.3}),
it follows that the kernel $K$ in \eqref{eq:euc-kern} is also positive definite.
Accordingly, the kernel
$Q(x,y) := (1 - \la x,y \ra)^{-s/2}$ on $\bS^n$
is positive definite (cf.\ Remark~\ref{rem:4.2}(b)).

Let $\mu_\bS$ be the $\OO_{n+1}(\R)$-invariant
measure on $\bS^n$, which, in stereographic coordinates, is given by
\[ d\mu_\bS(x) = \frac{2^n}{(1 + \|x\|^2)^n} \, dx. \]

\begin{rem} \mlabel{rem:hmu} Let $X$ be a locally compact space. A Radon measure 
$\mu$ on $X \times X$ is called {\it positive definite} if 
\[ \int_{X \times X} f(x) \oline{f(y)}\, d\mu(x,y) \geq 0 \quad \mbox{ for } \quad 
f \in C_c(X).\] 
Completing $C_c(X)$ with respect to the scalar product 
\[ \la f, g \ra_\mu := 
\int_{X \times X} f(x) \oline{g(y)}\, d\mu(x,y) \] 
then leads to a Hilbert space $\cH_\mu$ (cf.\ Remark~\ref{rem:2.2}). 
If $X$ is compact, then $\cH_\mu$ can most naturally 
be realized as a subspace 
of the space $\cM(X) := C(X)'$ of Radon measures on $X$ via the inclusion 
\[ T_\mu \: \cH_\mu \into \cM(X), \quad 
T_\mu(f)(g) = \la f, g \ra_\mu, \quad 
T_\mu(f)(E) = \int_{X \times E} f(x)\, d\mu(x,y).\] 

In the special case where 
$\mu$ has a density $\rho$ with respect to a product measure 
$\mu_X \otimes \mu_X$, the measure $T_\mu(f)$ can be written as 
$\Gamma_\mu(f)\mu_X$ with 
\[ \Gamma_\mu(f)(y) = \int_X f(x) \rho(x,y)\, d\mu_X(x).\] 
\end{rem}

\begin{lem}\label{le-Qintegrble}
For $Q(x,y) := (1 - \la x,y\ra)^{-s/2}$ and 
$0 \leq s < n$, the measure
\begin{equation} \label{eq:meas}
d\mu(x,y) := Q(x,y)\, d\mu_\bS(x) d\mu_\bS(y)
\end{equation}
on $\bS^n \times \bS^n$ is a finite positive Radon measure
which is positive definite. 
\end{lem}

\begin{prf} First we observe that the
$\OO_{n+1}(\R)$-invariance of the
measure $\mu_\bS$ and the kernel~$Q$
implies that the function $F(y) := \int_{\bS^n} Q(x,y)\, d\mu_\bS(x)$
is constant. For $y=e_0$, we find with \cite[Prop.~7.3.11]{vD09}
and some $c > 0$:
\begin{eqnarray}  \label{eq:F(e_0)}
F(e_0) &=& \int_{\bS^n}Q(x,e_0)\, d\mu_\bS(x) =c\int_{-1}^1 (1-t)^{-s/2}(1-t^2)^{(n-2)/2}\, dt\notag\\
& =&c\int_{-1}^1(1+t)^{(n-2)/2}(1-t)^{(n-2-s)/2}\, dt
\end{eqnarray}
and this integral is finite if and only if $s<n$.
This implies that
\[ \mu(\bS^n \times \bS^n) = \int_{\bS^n} F(y)\, d\mu_\bS(y) < \infty.\]

To see that $\mu$ is positive definite, we
use the conformal map $\eta \: \R^n\to \bS^n$
and Remark~\ref{rem:4.2}(b),(c) to obtain 
\begin{align*}
&\int_{\bS^n \times \bS^n} f(x)\oline{f(y)} \, d\mu(x,y) 
= c\int_{\R^n\times\R^n} \frac{f(\eta(x))}{(1+\|x\|^2)^{n+s/2}}
\frac{\oline{f(\eta(y))}}{(1+\|y\|^2)^{n+s/2}}\, \frac{dxdy}{\|x-y\|^s}\ge 0,
\end{align*}
where we have used that the
the kernel $\|x-y\|^{-s}$ is positive definite for~$0<s<n$
(Lemma~\ref{lem:1.3}).
\end{prf}

\begin{lem} For $f\in C(\bS^n)$, the function
\[ \Gamma(f) \: \bS^n \to \C, \quad \Gamma(f)(y) :=
 \int_{\bS^n}f(x)Q(x,y)\, d\mu_\bS(x)\]
is continuous.
\end{lem}

\begin{prf} It is enough to show that
$F(g):=\int_{\bS^n}f(x)Q(x,g.e_0)\, d\mu_\bS(x)$
is a continuous function on $\OO_{n+1}(\R )$. In view
of $Q(x,g.y)=Q(g^{-1}.x,y)$, we have
\[F(g)=\int_{\bS^n}f(g.x)(1-\la x, e_0\ra )^{-s/2}\, d\mu_\bS(x)\, .\]
The claim now follows
by $|f(g.x)(1-\la x, e_o\ra )^{-s/2}|\le \|f\|_\infty (1-\la x,e_o \ra)^{-s/2}$  and this function is integrable by \eqref{eq:F(e_0)} and $g\mapsto |f(g.x)|$ is continuous.
\end{prf}

%

\begin{defn} (The Hilbert spaces $\cH_s$, $0 < s < n$) 
Lemma~\ref{le-Qintegrble}
implies that $\mu$ is positive definite, 
so that we obtain with Remark~\ref{rem:hmu} a corresponding Hilbert space
$\cH_s := \cH_\mu \subeq \cM(\bS^n)$ of measures. 
It is the completion of the
range of the map
\[ T_\mu  \: C(\bS^n) \to  \cM(\bS^n), \quad
f \mapsto \Gamma(f) \mu_\bS,\]
with respect to the inner product
\begin{equation}
  \label{eq:inprogamma}
\la \Gamma(f_1) \mu_\bS, \Gamma(f_2) \mu_\bS \ra = \la f_1, f_2 \ra_\mu.
\end{equation}
\end{defn}

\begin{lem} \mlabel{lem:4.4b}
The prescription
\[ \pi_s(g)\nu
:= J_{g^{-1}}^{s/2-n} \cdot g_*\nu\]
defines a unitary representation of $G = \OO_{1,n+1}^+(\R)$
on the Hilbert subspace
$\cH_s \subeq \cM(\bS^n)$. 
\end{lem}

\begin{prf} 
For $g \in G$ we recall from Remark~\ref{rem:4.2}(a) that 
\begin{equation}
  \label{eq:trafoq}
Q(g.x,g.y) = J_g^{-s/2}(x) J_g^{-s/2}(y) Q(x,y). 
\end{equation}
The  measure $\mu_\bS$ transforms according to
\begin{equation}
  \label{eq:trafomeas2}
g_* \mu_\bS = J_{g^{-1}}^n \cdot \mu_\bS \quad \mbox{ for } \quad g \in G.
\end{equation}
Combining \eqref{eq:trafoq} with \eqref{eq:trafomeas2}
leads to
\begin{equation}
  \label{eq:trafomeas}
d(g_*\mu)(x,y) = J_{g^{-1}}(x)^{n-s/2} J_{g^{-1}}(y)^{n-s/2} d\mu(x,y).
\end{equation}
Therefore the action of $G$ on $C(\bS^n)$ defined by
\begin{equation}
  \label{eq:actfunc}
g.f := J_{g^{-1}}^{n - s/2} \cdot g_*f = g_*(J_g^{s/2-n} \cdot f) 
\end{equation}
satisfies
\begin{align*}
&\ \ \ \int_{\bS^n \times \bS^n} (g.f_1)(x) \cdot \oline{g.f_2(y)}\, d\mu(x,y)\\
&=  \int_{\bS^n \times \bS^n}
 (J_g^{s/2-n}f_1)(x) (J_g^{s/2-n}f_2)(y)\, d(g^{-1}_*\mu)(x,y) \\
&=  \int_{\bS^n \times \bS^n}  f_1(x)\oline{f_2(x)}\, d\mu(x,y),
\end{align*}
i.e., the sesquilinear form defined by the measure $\mu$
on $C(\bS^n)$ is invariant under the $G$-action. This in turn
implies that $T_\mu$ is equivariant with respect to the
action on $\cH_s \subeq \cM(\bS^n)$ given by
\begin{align*}
(g.\nu)(f)
&:= \nu(g^{-1}.f)
= \nu\big((g^{-1})_*(J_{g^{-1}}^{s/2-n} \cdot f)\big)
= (g_*\nu)(J_{g^{-1}}^{s/2-n} \cdot f)
= (J_{g^{-1}}^{s/2-n} \cdot g_*\nu)(f).\qedhere
\end{align*}
\end{prf}

\begin{rem} The preceding proof shows that the map
\[ T_\mu \: C(\bS^n) \to \cH_s, \quad f \mapsto \Gamma(f) \mu_\bS \] 
is $G$-equivariant with respect to the action on $C(\bS^n)$ given by 
\eqref{eq:actfunc}: $g.f := J_{g^{-1}}^{n - s/2} \cdot g_*f.$ 
As $\pi_s(g)(\Gamma(f)\mu_\bS) = J_{g^{-1}}^{s/2}\cdot g_*\Gamma(f)\mu_\bS$, it follows that 
\begin{equation}\label{eq:gammaeq}
\Gamma(g.f) = \Gamma(J_{g^{-1}}^{n - s/2} \cdot g_*f) 
= J_{g^{-1}}^{s/2} \cdot g_*\Gamma(f). 
\end{equation}
\end{rem}

\begin{rem} For $n = 1$, every diffeomorphism of $\bS^1$ is conformal,
so that one may expect that the representation $\pi_s$ on $\cH_s$
can be extended from
$\OO_{1,2}(\R)$ to the group of all diffeomorphisms of $\bS^1$.
Clearly, the conformal cocycle $J_g(x) := |\dd g(x)|$ is well-defined,
but we also need the projective invariance of the kernel $Q$, i.e., 
the invariance up to a cocycle.
In view of Remark~\ref{rem:4.2}(d), the projective invariance of $Q$
under some $g \in \Diff(\bS^1)$
implies the preservation of the absolute value of the cross ratio.
This implies that $g$ is fractional linear. In fact, composing
with a suitable fractional linear map, we may assume that
$g$ fixes $0$, $1$ and $\infty$. Then the invariance of the absolute
value of the cross ratio leads to
$|g(x)| = |x|$ for $x \in \R$ (in stereographic coordinates).
Now $g(1) = 1$ implies that $g = \id_{\bS^1}$.
\end{rem}

The following abstract lemma is useful to identify smooth and
distribution vectors for the representation $(\pi_s, \cH_s)$.

\begin{lem} \mlabel{lem:contlem}
Let $X$ be a topological vector space and
$F_1 \subeq F_2 \subeq X$ be two linear subspaces, both
carrying Fr\'echet topologies for which the inclusions
$F_j \to X$, $j =1,2$, are continuous. Then the
inclusion $\iota \: F_1 \to F_2$ is also continuous.
\end{lem}

\begin{prf} Since $F_1$ and $F_2$ are Fr\'echet spaces, it suffices
to verify that the graph of $\iota$ is closed.
So let $(v_n, \iota(v_n)) \to (v,w)$ in $F_1 \times F_2$.
Then $v_n \to v$ in $F_1$ implies $v_n \to v$ in $X$.
We also have $v_n \to w$ because the inclusion $F_2 \to X$
is continuous, and this leads to $w= v$. Hence the graph of $\iota$ is closed.
\end{prf}

\begin{lem} For each $f \in C^\infty(\bS^n)$, the element 
$\Gamma(f)\mu_\bS \in \cH_s$ is a smooth vector. Conversely, 
every smooth vector $\nu \in \cH_s^\infty$ is of the form 
$\nu_f = f \mu_\bS$ for a smooth function $f \in C^\infty(\bS^n)$ 
and the so obtained map 
$\cH_s^\infty \to C^\infty(\bS^n), \nu_f \mapsto f$ is continuous. 
In particular, its adjoint defines a linear map 
\[ \Psi \: \cD'(\bS^n) \to \cH_s^{-\infty}, \quad 
\Psi(D)(\nu_f) := D(f).\] 
\end{lem}

\begin{prf} In view of 
 $\|\Gamma(f)\mu_\bS\|^2 = \la f, f \ra_\mu$ (cf.\ \eqref{eq:inprogamma}), 
the map 
\[  C(\bS^n) \to \cH_s, \quad f \mapsto \Gamma(f) \mu_\bS \] 
is continuous. We have seen in \eqref{eq:gammaeq} 
that it is $G$-equivariant with respect to the action on 
$C(\bS^n)$ by $g.f = J_{g^{-1}}^{n-s/2} g_*f$. As the cocycle 
$J$ is smooth, smooth functions $f$ on $\bS^n$ have smooth orbit 
maps for this action, so that $\Gamma(f)\mu_\bS \in \cH_s^\infty$.

{}From \cite[Thm.~3.3]{DM78} we know that
$\cH_s^\infty$ is spanned by $\pi_s(\cD(G))\cH_s$.
Accordingly, we obtain for the maximal compact subgroup
$K \cong \OO_{n+1}(\R)$ of $G = \OO_{1,n+1}^+(\R)$ that the corresponding space
$\cH_s^\infty(K)$ of smooth vectors for $K$
is spanned by $\pi_s(\cD(K))\cH_s$. Since $K$ acts on $\bS^n$
by isometries, $J_k = 1$ for every $k \in K$, so that
$\pi_s(k)\nu = k_* \nu$.
As $\cH_s$ is realized in $\cM(\bS^n)$, which can be identified
with a subspace of $\cM(K)$ because $K$ acts transitively on $\bS^n$,
it follows that
\[ \cH_s^\infty \subeq  \Spann(\pi_s(\cD(K))\cH_s)
\subeq \cD(K) * \cM(\bS^n) \subeq C^\infty(\bS^n)\cdot \mu_\bS \]
(cf.\ \cite[Prop.~A.2.4.1]{Wa72}). From the continuity of the inclusions
$\cH_s^\infty \to \cH_s^\infty(K) \to\cM(\bS^n)$ and Lemma~\ref{lem:contlem},
it now follows that the linear map
$\cH_s^\infty \into C^\infty(\bS^n)$ 
is also continuous. Its adjoint therefore defines a continuous linear map
$\cD'(\bS^n) \to \cH_s^{-\infty}.$
\end{prf}

\begin{rem} \mlabel{rem:5.10} 
The representations $(\pi_s,\cH_s)_{0 < s < n}$ are complementary series 
representations. Those representations are well known and arise from intertwining
operators, usually given by singular integral operators, \cite{KS71}. The reader finds  a
detailed discussion in Section 7.5 of \cite{vD09}.
\end{rem}

\section{Reflection positivity for $\pi_s$} \mlabel{sec:5}

In this section we study reflection positivity for the 
complementary series representations 
$(\pi_s, \cH_s)_{0 < s < n}$ introduced in the preceding section. 
It turns out that these representations of the conformal group 
$\OO_{1,n+1}(\R)$ on $\R^n$ provide a natural context for relating 
reflection positivity, resp., 
positive definiteness of kernels of the form 
$K(x,\tau y)$ (cf.\ Example~\ref{ex:ker-refpos}) on open half spaces 
and open balls. 

\subsection{The half space picture}
\mlabel{subsec:6.1}

Formulas for kernels and questions of positive definiteness
turn out to be simpler for the open half space $\R^n_+$ compared to the open unit 
ball~$\cD$. In particular, this connects much better with involutive
semigroups and integral representations of positive definite functions
thereon. We therefore discuss this case first.

As before, let
\[ \Omega = \{ z \in \R^{n} \: z_0 > 0, [z,z] > 0\} \subeq \R^{n}\]
be the forward light cone and
$T_\Omega := \Omega + i \R^{n} \subeq \C^{n}$ be the corresponding
tube domain. Then it is easy to show that the complex bilinear
extension of $[\cdot,\cdot]$ to $\C^{n}$ satisfies
\[ \Delta(z) := [z,z] \not= 0 \quad \mbox{ for } \quad z \in T_\Omega.\]
Since $T_\Omega$ is simply connected, there exists a
holomorphic function $\log [z,z]$ taking the value
$0$ in $e_0$. We thus obtain on
this domain for each $\lambda \in \R$ a holomorphic function
\[ \Delta^{-\lambda} \: T_\Omega \to \C^\times, \quad
z \mapsto e^{-\lambda \log [z,z]}. \]

Since $(x + iy)^* := x - iy$ defines on $T_\Omega$ the structure
of a complex involutive semigroup, the function $\Delta^\lambda$
defines a kernel $(z,w) \mapsto \Delta^{-\lambda}(z + w^*)$
on $T_\Omega$. From \cite[Thm.~XIII.2.7]{FK94},
we know that this kernel is positive definite, i.e.,
$\Delta^{-\lambda}$ is a positive definite function on
$T_\Omega$, if and only if
\begin{equation}
  \label{eq:wallset}
\lambda \in
\cW := \Big\{0,\frac{n-2}{2}\Big\} \cup \Big]\frac{n-2}{2}, \infty\Big[.
\end{equation}
{}From \cite[Lemma~X.4.4 and p.~262]{FK94} it now follows
that this is equivalent to the positive definiteness of the 
function $\Delta^{-\lambda}$ on the light cone $\Omega$. 

The domain
\[ \Omega^c
= \{ z \in T_\Omega \: z_0 \in \R, z_1,\ldots, z_{n-1} \in i \R\}
= \R_+ \times i \R^n\]
is a totally real submanifold and
also an involutive subsemigroup of $T_\Omega$.
Therefore Theorem~\ref{thm:a.2.1} implies that the restriction
of $\Delta^{-\lambda}$ to $\Omega^c$ is a positive definite function
if and only if \eqref{eq:wallset} holds.

Define $\iota \: \R^{n} \to \C^n, x \mapsto (x_0, ix_1,\ldots, ix_{n-1})$.
Then $\iota$ restricts to an isomorphism of involutive semigroups
$\R^n_+ \to \Omega^c$, where the involution on $\R^n_+$ is given by
\[ (x_0, x_1, \ldots, x_{n-1})^\sharp = (x_0, -x_1,\ldots, -x_{n-1}). \]
In view of $\Delta(\iota(x)) = [\iota(x),\iota(x)] = \|x\|^2,$
we have for $\lambda = s/2$ the relation
$(\iota^*\Delta^{-\lambda})(x) = \|x\|^{-s}.$
This leads to:

\begin{prop} \mlabel{prop:wallset}
The function $\|x\|^{-s}$ is positive definite on the involutive semigroup 
$(\R^n_+,\sharp)$ if and only if $s = 0$ or $s \geq \max\{0,n-2\}$.
\end{prop}

\begin{prf} For $n \geq 2$ this follows from the preceding
discussion, and for
$n = 1$ we have already
seen in Example~\ref{ex:2.3} that the function $x^{-s}$ is positive
definite on $\R_+$ for\break $s \geq 0$.
\end{prf}

\subsection{The ball picture}

\begin{prop} \mlabel{prop:2.5b} The kernel
$$ R(x,y) := (1 - 2 \la x, y \ra + \|x\|^2 \|y\|^2)^{-s/2} $$
on the open unit ball $\cD \subeq \R^n$ is positive definite if and only if
\[ s = 0 \quad \mbox{ or } \quad s \geq \max(0,n - 2).\]
\end{prop}

\begin{prf} Open half spaces and open balls
in $\R^{n}$ are conformally equivalent. The equivalence is obtained
by rotation a lower hemisphere (open ball) in the conformal compactification
into a right hemisphere, which corresponds to a half space.
Let $c \in \OO_{n+1}(\R)$
be an involution exchanging $\R^n_+$ with the unit ball $\cD$
and $\tilde c$ be the corresponding map on $\R^n$ defined by
$\eta \circ \tilde c = c \circ \eta$ (Remark~\ref{rem:4.3}).
Then $J_c = 1$, so that $c$ leaves the
kernel $Q$ invariant.

Let $\tau \in G$ be the element inducing on $\R^n$ the
reflection $\tilde\tau$ in $\partial\R^n_+ \cong \R^{n-1}$,
and $\sigma$ be the reflection in $\partial\cD$. Then
$\tau = c\sigma c$. Now the invariance of the kernel $Q$ on $\bS^n$ 
under $c$ leads to 
\begin{align*}
K(\tilde\tau(x),y)
&= Q(\eta(\tilde\tau(x),\eta(y)) = Q(\tau\eta(x),\eta(y))
= Q(\sigma c\eta(x),c\eta(y)) \\
&= Q(\eta(\tilde\sigma\tilde c(x)),\eta(\tilde c(y)))
= K(\tilde\sigma(\tilde c(x)), \tilde c(y)).
\end{align*}
Clearly, this kernel is positive definite on $\R^n_+$ if and only
if $K(\tilde\sigma(x),y)$ is positive definite on $\cD$. 

{}From $\tilde\sigma(x) = \frac{x}{\|x\|^2}$ we obtain 
$$ \|x\|^2 \|\tilde\sigma(x)- y\|^2
=   \Big\|\frac{x}{\|x\|} - \|x\| y\Big\|^2
=  1 - 2 \la x, y \ra + \|x\|^2 \|y\|^2.   $$
Therefore the pullback of the kernel $Q(\sigma(x),y)$ under
$\eta$ is given by 
\begin{eqnarray} \label{eq:twistker}
K(\tilde\sigma(x),y)
&=& (1 + \|\sigma(x)\|^2)^{s/2} \|\sigma(x) -y\|^{-s}
(1 + \|y\|^2)^{s/2} \notag \\
&=& (1 + \|x\|^2)^{s/2} \|x\|^{-s} \|\sigma(x) -y\|^{-s}
(1 + \|y\|^2)^{s/2}  \\
&=& (1 + \|x\|^2)^{s/2}
{(1 - 2 \la x, y \ra + \|x\|^2 \|y\|^2)^{-s/2}}
(1 + \|y\|^2)^{s/2}.\notag
\end{eqnarray}
We conclude that the kernel $K(\tilde\tau(x),y)$ 
is positive definite if and only if $R$ is positive definite. 

Next we observe that 
\begin{align*}
K(\tilde\tau(x),y)
&= (1 + \|\tilde\tau(x)\|^2)^{s/2}
\|\tilde\tau(x) -y\|^{-s} (1 + \|y\|^2)^{s/2} \\
&= (1 + \|x\|^2)^{s/2} \big((x_0+y_0)^2 + \|x'-y'\|^2\big)^{-s/2}
(1 + \|y\|^2)^{s/2}
\end{align*}
is positive definite if and only if the kernel
\[ (x,y) \mapsto \big((x_0+y_0)^2 + \|x'-y'\|^2\big)^{-s/2} \]
is positive definite on $\R^n_+$. This means that the function
$\Delta(x) = \|x\|^{-s}$
is positive definite on the involutive semigroup $(\R_n^+,-\tau)$.
In view of Proposition~\ref{prop:wallset}, this is equivalent to
$s = 0$ or $s\geq \max(0,n-2)$.
\end{prf}

\begin{rem} For an alternative proof of the preceding proposition, 
one can observe that the unit ball $\cD \subeq \R^n$ is a totally 
real submanifold of the Lie ball $\cD_\C \subeq \C^n$. Since $R$ is a power 
of the Bergman kernel on the Lie ball which is biholomorphic 
to the tube domain $T_\Omega$, one can also derive our result from 
\cite[Thm.~XIII.2.7]{FK94} by using Theorem~\ref{thm:a.2.1} and 
$\cD = \cD_\C \cap \R^n$. 
\end{rem}

\begin{ex} (a) For $n=1$ we have $R(x,y) = (1 - 2 xy + x^2 y^2)^{-s/2}
= (1 - xy)^{-s}$, which is positive definite for each $s \geq 0$.

(b) For $n = 2$ the kernel $R$ is also positive definite for $s \geq 0$. This is basically 
due to the fact that it corresponds to the positive definiteness of the function 
$\Delta^{-s/2}(x) = [x,x]^{-s/2}$ on the open light cone $\Omega$ (cf.\ Subsection~\ref{subsec:6.1}). 
>From 
\[ [x,x] = x_0^2 - x_1^2 = (x_0 - x_1)(x_0 + x_1) \]  we obtain a factorization 
of this function, so that the positive definiteness of the functions 
$(x_0 \pm x_1)^{-s/2}$ on $\Omega$ for $s \geq 0$ implies that $R$ is positive definite. 
\end{ex}

\subsection{Reflection positivity on $\cH_s$}
Let us fix some notation used in this section. We write $G=\OO^+_{1,n+1}(\R)$ and
$\g=\fo_{1,n+1}(\R)$ for the Lie algebra of $G$.  Let $\cDt\subeq \R^n$ 
be the interior of the unit ball and $\cD=\eta(\cDt)$ be its image in $\bS^n$. 
We consider the corresponding compression semigroup 
\[ S_\cD := \{ g \in G \:  g\cD \subeq \cD\}  \]
in $G$. Let $\sigma \in G$ be the element implementing 
the conformal reflection $\tilde{\sigma} (x)= x/\|x\|^2$ in
$\bS^{n-1}=\partial \cDt$. On $\bS^n$, it acts by 
\[ (x_0,x_1,\ldots ,x_n)\mapsto (-x_0,x_1,\ldots ,x_n).\]
According to (\ref{eq:confact}), it is realized by the diagonal matrix
\[ \sigma := \diag(1,-1,1,\ldots ,1) \in \OO_{1,n+1}(\R).\]
The corresponding involution on $\OO_{1,n+1} (\R )$ is given by the
conjugation by the same matrix.   $\Ad (\sigma)$ defines 
an involution on $\g$, so that $\g$ decomposes into
$\pm 1$-eigenspaces $\fg = \fh \oplus \fq$, where  $\fh = \Fix(\Ad(\sigma))$
and $\fq = \Fix(-\Ad(\sigma))$. The Lie algebra $\fh$ is isomorphic 
to $\fo_{1,n}(\R)$. 

\begin{prop} \mlabel{prop:compsem} Let
$\sigma \in G$ be the conformal reflection in $\partial\cD$. Then
\[ S_\cD = H \exp(C),\quad \mbox{ where } \quad
H = \OO_{1,n}^+(\R)\subeq G^\sigma,\]
and $C \subeq \fq$ is a closed convex $\Ad(H)$-invariant cone.
\end{prop}

\begin{prf} Let $\tilde G_0$ be
the simply connected covering group of the identity component
$G_0 = \OO_{1,n+1}(\R)_0 = \SO_{1,n+1}^+(\R)$
and $\hat\sigma$ denote the lift of the involution defined by
conjugation with $\sigma$ on $G_0$ to its universal covering
group~$\tilde G_0$.
Then the group $(\tilde G_0)^{\hat\sigma}$ is connected
(\cite[Thm.~IV.3.4]{Lo69}). From
\cite[Thm.~VI.11, Rem.~VI.12]{HN95} we thus derive
that
\[ S_\cD \cap G_0 = \SO_{1,n}^+(\R)\exp(C)\]
for a proper closed convex cone $C \subeq \fq$ which
is invariant under $\Ad(\SO_{1,n}^+(\R))$.

Since the group $G$ has two connected components
determined by the values of the determinant,
and
\[ H = \OO_{1,n}^+(\R) \subeq G^\sigma
\cong \OO_1(\R) \times \OO_{1,n}^+(\R)
= \{\1,\sigma\} \times \OO_{1,n}^+(\R)\]
likewise does, it follows that $G = H G_0$. Clearly,
$G^\sigma$ preserves the fixed point set $\partial \cD$ of $\sigma$
in $\bS^n$ and $\sigma$ exchanges the two connected components of
its complement, which are both preserved by $H$. 
This shows that $H \subeq S_\cD$.

Pick $h_0 \in H$ with $\det(h_0) = -1$.
If $s \in S_\cD \setminus G_0$, then
$h_0^{-1}s \in S_\cD \cap G_0$, so that
$S_\cD = H (S_\cD \cap G_0) = H \exp(C).$
\end{prf}

Let $K=\OO_n(\R)\subset \OO^+_{1,n+1}(\R)$ and $H_K :=H\cap K\simeq \OO_{n-1}(\R)$.
On the subgroup $A \cong (\R_+,\cdot) \subeq G$, acting on
$\R^n$ by multiplication with positive scalars (cf.\ \eqref{eq:affconf}),
the involution
$\sigma$ acts by inversion, so that there exists a Lie algebra element $x_o$ such that
$A \cap S = \exp(\R_+ x_o)$
is a one-parameter semigroup acting
as $(]0,1[,\cdot)$ on $\R^n$ which commutes with $H_K$.
Moreover, we have $C=\Ad (H)\R_+ x_o$, (cf.~\cite[Sect.~4.4]{HO96}).

On the Hilbert space $\cH_s \subeq \cM(\bS^n)$, we consider the
involution $\theta  := \pi_s(\sigma)$, where
$\sigma \in G = \OO_{1,n+1}^+(\R)$ is as above.

\begin{lem} \mlabel{lem:hs-repo} 
If $s = 0$ or $s \geq n -2$, then the closed subspace 
\[ \tilde\cE_+ := \oline{\{\Gamma(f)\mu_\bS \: \supp(f) \subeq \cD\}} \subeq \cH_s \] 
is $\theta $-positive. 
\end{lem}

\begin{prf} As $\sigma$ acts isometrically on $\bS^n$, we have 
$J_\sigma = 1$, so that 
\[ \theta  \nu = \pi_s(\sigma)\nu = \sigma_* \nu.\]
For $f \in C_c(\cD)$  we thus find 
\begin{align*}
\la \theta  \Gamma(f) \mu_\bS, \Gamma(f) \mu_\bS \ra 
&=\la \Gamma(\sigma_* f) \mu_\bS, \Gamma(f) \mu_\bS \ra \\
&= \int_{\bS^n} \int_{\bS^n}
f(\sigma(x)) \oline{f(y)} Q(x,y) \, d\mu_\bS(x)d\mu_\bS(y) \\
&= \int_{\cD} \int_{\cD} f(x) \oline{f(y)} Q(\sigma(x),y)
\, d\mu_\bS(x)d\mu_\bS(y).
\end{align*}
To evaluate this integral in stereographic coordinates, 
we recall from \eqref{eq:twistker} the proof of Proposition~\ref{prop:2.5b} 
that  the pullback of the kernel $Q(\sigma(x),y)$ under
$\eta$ is 
\[ (1 + \|x\|^2)^{s/2}
{(1 - 2 \la x, y \ra + \|x\|^2 \|y\|^2)^{-s/2}}
(1 + \|y\|^2)^{s/2}.\] 
We thus obtain with a positive constant $c'$:
\[ \la \theta \nu, \nu \ra = c'
 \int_{\cD} \int_{\cD}
\frac{f(x) \oline{f(y)} (1 + \|x\|^2)^{s/2}(1 + \|y\|^2)^{s/2}}
{(1 - 2 \la x, y \ra + \|x\|^2 \|y\|^2)^{s/2}}
\, \frac{dx}{(1 + \|x\|^2)^n} \frac{dy}{(1 + \|y\|^2)^n}.\]
That this expression is non-negative for
$f \in C_c(\cD)$ follows from Proposition~\ref{prop:2.5b}.
\end{prf}

\begin{thm} \mlabel{thm:6.7} 
For $x\in \bS^{n-1}=(\bS^n)^{\sigma}$ let $\delta_x\in \cH_s^{-\infty}$ be
the delta measure in $x$. The triple 
$(\pi_s,\cH_s,\delta_x)$ is a reflection positive distribution cyclic representation
for $(G, \tau, S_\cD^0)$ if $s=0$ or $n-2\le s< n$.
\end{thm}

\begin{prf} If
$\phi\in \cD (G)$ then $\pi_s^{-\infty}(\phi)\delta_x\in \cH^\infty_s\subeq 
\cD (\bS^n)\mu_\bS$. Thus $\pi_s^{-\infty}(\phi )\delta_x =\phi^\flat\mu_\bS$ for
some uniquely determined $\phi^\flat\in \cD (\bS^n)$. 

To determine this function, we first have to identify 
the distribution vector $\delta_x$ as a measure on $\bS^n$. 
Since $\la \delta_x, \Gamma(f) \mu_\bS \ra = \oline{\Gamma(f)(x)}$ 
for $f \in C^\infty(\bS^n)$, the measure corresponding to $\delta_x$ 
is
\[ \Gamma^*(\delta_x) = Q_x \cdot \mu_\bS\]  
and therefore $\pi_s^{-\infty}(\phi)\delta_x$ corresponds to the measure 
\[ \int_G \phi(g) \pi_s(g)(Q_x \cdot \mu_\bS)\, d\mu_G(g) 
= \left( \int_G \phi(g) J_{g^{-1}}^{s/2} g_* Q_x \, d\mu_G(g)\right) \cdot \mu_\bS,\] 
which means that 
\[ \phi^\flat  =  \int_G \phi(g) J_{g^{-1}}^{s/2} g_* Q_x \, d\mu_G(g).\] 

Let $G_x$ be the stabilizer of $x$ in $G$.  Then
$G/G_x \simeq \bS^n$ via $gG_x\mapsto g.x$.  The quasi-invariant measure
$\mu_{S}$ on $\bS^n$ and the left invariant measure on $G$ are related by\begin{footnote}{This formula is most easily verified by showing that 
$f \mapsto \int_{\bS^n} \int_{G_x}f(kp)J_{kp}^{-n}(x)\, d\mu_{G_x}( p )\, d\mu_{S}(k.x)$ 
defines a left invariant integral on $G$.}
\end{footnote}
\[\int_G f(g)J_g^n (x)\, d\mu_G(g)=\int_{\bS^n} \int_{G_x}f(kp)\, d\mu_{G_x}( p )\, d\mu_{S}(k.x)\, .\]
{}From \eqref{eq:trafoq} we have $Q(g.x,g.y) = J_g^{-s/2}(x) J_g^{-s/2}(y) Q(x,y)$ and hence
\begin{align*}
(g_*Q_x)(y)  
&=  Q(g^{-1}.y, x) 
= Q(y,g.x) J_g^{s/2}(x)J_g^{s/2}(g^{-1}.y) \\
& =  Q(y,g.x) J_g^{s/2}(x)J_{g^{-1}}^{-s/2}(y).
\end{align*}
This leads to 
\begin{align*}
\phi^\flat
&=  \int_G \phi(g) J_g^{s/2}(x) Q_{g.x}\,  d\mu_G(g)
=  \int_{\bS^n} \int_{G_x} J_{gp}^{s/2-n} ( x) \phi(gp)\, d\mu_{G_x} (p ) Q_{g.x}\,  d\mu_\bS(g.x)\, .
\end{align*}
This means that 
\begin{equation}\label{eq:6.40}
\phi^\flat =  \Gamma(\Phi) \quad \mbox{ with } \quad
\Phi(g.x)  = \int_{G_x} J_{gp}^{s/2-n}(x) \phi(gp )\, d\mu_{G_x} ( p )\, .
\end{equation}

As $J_g(x)>0$ for all $g$ it follows that the map $\cD (G)\to \cD (\bS^n)$, $\phi \mapsto \Phi$, is surjective and
hence that 
$\pi_s^{-\infty}(\cD (G))\delta_x$ is dense in $\cH_s$. Thus $\delta_x$ is a cyclic distribution vector. 

Suppose that $\supp (\phi )\subset S_\cD^0$ and assume that $\Phi (g.x)\not= 0$. Then there exists
$p\in G_x$ such that $gp\in \supp(\phi) \subeq S_\cD^0$. Hence $g.x\in S_\cD^0.x \subset \cD$. Thus 
$\supp \Phi \subset \cD$, and hence 
\begin{equation}\label{eq:6.41}
\pi^{-\infty}_s(\cD (S))\delta_x\subset \tilde\cE_+ 
\end{equation}
(cf.\ Lemma~\ref{lem:hs-repo}).  The claim now follows from Lemma \ref{lem:hs-repo}.
\end{prf}

\begin{rem} As $\OO_n(\R)$ acts transitively on $\bS^n$ we can in (\ref{eq:6.40}) assume that $g\in \OO_n(\R)$. Then $J_g(x)=1$ and
$J_{gp}(x)=J_g(p.x)J_p(x)=J_p(x)$. Thus
\[\Phi (g.x)=\int_{G_x} J_p(x)^{s/2-n}(x)\phi (g p) \, d\mu_{G_x}(p )\]
for $g\in \OO_n(\R)$.
\end{rem}

\begin{rem} If $\alpha$ is a positive linear combination, or even an integral with respect to a positive Borel measure,
of $\delta$-distributions supported on $\bS^{n-1}$, then \eqref{eq:6.41} shows that
$\pi_s^{-\infty}(\cD (S))\alpha \subset \tilde{\cE_+}$. 
According to Remark~\ref{rem:5.10} and \cite[p.119]{vD09}, 
the representation  $(\pi_s,\cH_s)$ is irreducible. Thus every nonzero distribution vector is cyclic. It then follows that $(\pi_s,\cH_s,\alpha)$
is a reflection positive distribution cyclic representation. In particular this holds for the measure $\mu_{\bS^{n-1}}$. We have
\[\Gamma^*(\mu_{\bS^{n-1}})  =\int_{\bS^{n-1}}Q_{x}\, d\mu_{\bS^{n-1}}(x)\, \mu_{\bS}\, .\]
\end{rem}

Let $\g=\h\oplus \fq$ be a symmetric Lie algebra corresponding to the involution $\tau$ on $G$ and 
$\g^c:=\h\oplus i\fq$ the $c$-dual symmetric Lie algebra. Let $G^c$ denote
the simply connected Lie group with Lie algebra $\g^c$. Then $\tau$ defines an involution
on $G^c$ and $(G^c)^\tau$ is connected (\cite[Thm.~IV.3.4]{Lo69} or \cite[Thm.~ 8.2,p. 320]{He78}). The symmetric space $G^c/(G^c)^\tau$ is the $c$-dual of $G/H$. 
As explained in \cite{JO00}, Sections~6 and 10, see also \cite{HN93}, \cite{JO98}, 
the reflection positivity and the L\"uscher--Mack 
Theorem now gives an irreducible highest weight (or positive energy) representation of the 
$c$-dual group $G^c$ on $\hat\cE$. 
Multiplying by $i$ in the
second coordinate transforms the Lorentz form
$[x,y]=x_0 y_0-x_1y_1-\ldots - x_{n+1}y_{n+1}$
into the form
\[ [x,y]_{2}=x_0y_0+x_1y_1-x_2y_2-\ldots -x_{n+1}y_{n+1}.\]
Hence the group $G^c$,  $c$-dual to $\OO^+_{1,n+1}(\R)$,
is locally isomorphic to $\OO_{2,n}(\R)$ . We point out that the condition that $s=0$ or $n-2\le s<n$ indicates that
this construction does not carry over to infinite dimension, or the duality between $\OO^+_{1,\infty}(\R)$ 
and $\OO_{2,\infty}(\R)$.
This is also reflected in the fact, that the group $\OO_{2,\infty}(\R)$ 
does not have any unitary highest weight representations  (\cite[pp.~276/277]{NO98}, 
\cite[Thm.~7.5]{Ne11b}).

\begin{rem} We have seen above that the open cones
$\Omega$ and $\R^n_+$, endowed with their natural semigroup involution,
can be viewed as real forms of the tube $T_\Omega = \Omega + i \R^n$,
endowed with its natural involution given by conjugation.
Every positive definite function $\phi \: \Omega \to \C$ extends
to $T_\Omega$ (\cite[Cor.\ to Thm.~4]{Sh84})
and hence restricts to $\R^n_+$ but the converse is
not true. For bounded positive definite functions
on $\Omega$,
we know that they are precisely the Laplace transforms
$\phi = \cL(\mu)$ of measures $\mu$ on the dual cone
$\hat\Omega = \Omega^\star$
whose Laplace transform is defined on $\Omega$.
Likewise bounded positive definite functions on $\R^n_+$ are
Fourier--Laplace transforms of measures on the closed half space
$\hat{\R^n_+} \supeq \Omega^\star$ and not all of them extend to
holomorphic positive definite functions on $T_\Omega$.
\end{rem}

\appendix

\section{Propagation of positive definiteness} \mlabel{app:b}

In this appendix we discuss some useful results providing
criteria for  kernels on complex manifolds
to be positive definite.

Let $M$ be a connected complex Fr\'echet
manifold, $\oline M$ its complex conjugate,
and $K \: M \times \oline M \to \C$ be a holomorphic function.
We call such functions {\it holomorphic kernels on $M$}.
A submanifold $\Sigma \subeq M$ is called {\it totally real}
if, for each point $s \in \Sigma$, there exists a holomorphic chart
$\phi \: U \to V_\C$, where $U$ is an open neighborhood of $s$ in $M$
and $V_\C$ is a complexification of the real locally convex space $V$,
such that $\phi(U \cap \Sigma) = \phi(U) \cap V$.

The following theorem generalizes \cite[\S 1.8]{Kr49}
from domains in $\C$ to complex Fr\'echet manifolds.

\begin{thm} \mlabel{thm:a.2.1}
For a holomorphic kernel $K$ on the connected complex Fr\'echet
manifold $M$ the following conditions are sufficient for
$K$ to be positive definite:
\begin{description}
  \item[\rm(i)] $K$ is positive definite on a non-empty open subset.
  \item[\rm(ii)] $K$ is positive definite on a non-empty totally real
submanifold.
\end{description}
\end{thm}

\begin{prf} (i) {\bf Step 1:} Let $\eset\not=U \subeq M$ be an open subset
on which $K$ is positive definite
and $\cH := \cH_{K\res_{U \times U}} \subeq \cO(U,\C)$
be the corresponding reproducing kernel Hilbert space.
We want to show
that the corresponding realization map
\[ \gamma \: U \to \cH, \quad \gamma(s) = K_s, \quad
K_s(t) = K(t,s) = \la K_s, K_t \ra \]
extends to an antiholomorphic map
$\gamma \: M \to \cH$.

If this is the case, then, for each $s \in U$, the map
$M \to \C, m \mapsto \la K_s, \gamma(m) \ra = \oline{\gamma(m)(s)}$
is holomorphic. For $m \in M$ we have
\[\oline{\gamma(m)(s)} = \oline{K_m(s)} = \oline{K(s,m)} = K(m,s),\]
so that the uniqueness of analytic continuation
implies that $\gamma(m)(s) = K_m(s)$, i.e.,
$\gamma(m) = K_m\res_U$. We conclude in particular that,
whenever $\gamma$ exists on some connected open subset $N \subeq M$
intersecting $U$, we necessarily have $\gamma(n) = K_n\res_{U}$
for each $n \in N$.

Next we note that the function
\[ M \times \oline M \to \C, \quad
(m,n) \mapsto \la \gamma(n), \gamma(m)\ra \]
is holomorphic and coincides for $(m,n) \in U \times U$ with
$K(m,n)$. By uniqueness of analytic continuation, we thus obtain
\[ K(m,n) = \la \gamma(n), \gamma(m)\ra \quad \mbox{ for } \quad
m,n \in M,\]
and hence that $K$ is positive definite on $N \times N$.

{\bf Step 2:} Suppose that $N_1$ and $N_2$ are two open connected
subsets of $M$ containing $U$ on which antiholomorphic
extensions
$\gamma^{N_1} \: N_1 \to \cH$ and
$\gamma^{N_2} \: N_2 \to \cH$ exist. For each $n \in N_1 \cap N_2$
Step $1$ implies that
$\gamma^{N^1}(n) = K_n\res_{U} = \gamma^{N_2}(n)$,
so that
$\gamma^{N_1} \res_{N_1 \cap N_2} = \gamma^{N_2} \res_{N_1 \cap N_2}$.
Therefore these two maps combine to a holomorphic map
\[ \gamma^{N_1 \cup N_2} \: N_1 \cup N_2 \to \cH.\] 

Let $N \subeq M$ be the union of all open connected subsets of
$M$ containing $U$ on which an antiholomorphic
extension of $\gamma$ exists. Then $\gamma$ extends to a holomorphic
map on $N$, and, in view of Step $1$, it only remains to show that $N = M$.

{\bf Step 3:} $N = M$. We argue by contradiction. Suppose
that $m_0 \in M$ is a boundary point of $N$. For each open
connected neighborhood $U$ of $m_0$ the intersection
$U \cap N$ is non-empty. Since $K$ is in particular continuous,
we may choose $U$ so small that $K$ is bounded on $U \times U$.
Fixing a local chart around $m_0$, we may further assume that
$U$ is biholomorphic to an open convex subset of a complex
locally convex space $X$. In the following we therefore consider
$U$ as such an open subset.

Then the arguments in the proof
of \cite[Thm.~5.1]{Ne11a} imply the existence of an open
$0$-neighborhood $U_1 \subeq X$ (only depending on the bound for $K$ on
$U \times U$) such that for every
point $m \in N \cap U$, the Taylor series of
$\gamma$ in $m$ converges in $m + U_1$ to a holomorphic
function which coincides with $\gamma$ in $(m + U_1) \cap N$.
Since $N$ intersects $m_0 - U_1$, there exists an
$m \in N$ with $m_0 \in m + U_1$. Therefore the holomorphic
function $\gamma\res_{(m + U_1)\cap N} \to \cH$ extends to a holomorphic
function $\hat\gamma \: m + U_1 \to \cH$. Step $2$ now implies that
$m + U_1 \subeq N$, contradicting $m_0 \in \partial N$.
This proves that $M = N$.

(ii) Let $\eset\not=\Sigma \subeq M$ be a totally real
submanifold on which $K^\Sigma := K\res_{\Sigma \times \Sigma}$
is positive definite. Then we obtain a reproducing
kernel Hilbert space
$\cH := \cH_{K^\Sigma}$ of functions on $\Sigma$. Again we want to show
that the corresponding realization map
\[ \gamma \: \Sigma \to \cH, \quad \gamma(s) = K_s^\Sigma, \quad
K_s^\Sigma(t) = K^\Sigma(t,s) = \la K_s^\Sigma, K_t^\Sigma \ra \]
extends to an antiholomorphic map
$\gamma \: M \to \cH$.

As in (i), any holomorphic extension $\gamma \: U \to \cH$ to
an open connected subset intersecting $\Sigma$
satisfies $\gamma(u)\res_\Sigma = K_u\res_{\Sigma}$ for
$u \in U$ and
\begin{equation}
  \label{eq:KU}
\la \gamma(u), \gamma(v) \ra = K(v,u)\quad \mbox{ for } \quad
u,v, \in U.
\end{equation}

Since $K$ is real analytic on $\Sigma \times \Sigma$, it follows from
\cite[Thm.~5.1]{Ne11a} that the map
$\gamma \: \Sigma \to \cH$ is analytic. By definition,
each point $s \in \Sigma$ has a connected neighborhood
to which $\gamma$ extends holomorphically. The preceding argument
shows that all these extensions can be patched together, so that
$\gamma$ extends holomorphically to an open neighborhood $U$ of
$\Sigma$. Now \eqref{eq:KU} implies that $K$ is positive
definite on $U$. Therefore (ii) follows from (i).
\end{prf}


\begin{thebibliography}{aaaaaaaa}
\bibitem[Ba78]{Ba78} Bauer, H., ``Wahrscheinlichkeitstheorie und
Grundz\"uge der Ma\ss{}theorie,'' Walter de Gruyter, Berlin, 1978

\bibitem[BG11]{BG11} Birth, L., and H. Gl\"ockner, {\it 
Continuity of convolution of test functions on Lie groups}, 
arXiv:math.FA/1112.4729v3 

\bibitem[Bl98]{Bl98}  Blackadar, B.,
``$K$-theory for Operator Algebras,'' 2nd Ed., Cambridge Univ. Press, 1998

\bibitem[Br93]{Br93} Bredon, G.\ E.,
``Topology and Geometry,'' Graduate Texts in
Mathematics {\bf 139}, Springer-Verlag, Berlin, 1993

\bibitem[vD09]{vD09} van Dijk, G., ``Introduction to Harmonic
Analysis and Generalized Gelfand Pairs,''
Studies in Math. {\bf 36} (2009), de Gruyter, Berlin, 2009

\bibitem[DM78]{DM78} Dixmier, J., and P. Malliavin,
{\it Factorisations de fonctions et de vecteurs ind\'efiniment
diff\'erentiables}, Bull. Soc. math., 2e s\'erie {\bf 102} (1978),
305--330

\bibitem[FK94]{FK94} Faraut, J., and A. Koranyi,
``Analysis on Symmetric Cones",
Oxform Mathematical Monographs, Oxford University Press, 1994

\bibitem[FL10]{FL10} Frank, R.~L., and E.~H. Lieb, {\it Inversion positivity and the 
sharp Hardy-Littlewood-Sobolev inequality}, Calc. Var. Partial Differential Equations {\bf 39:1-2} 
(2010), 85--99

\bibitem[FL11]{FL11} ---, {\it Spherical reflection positivity and the 
Hardy-Littlewood-Sobolev inequality}, in ``Concentration, functional inequalities 
and isoperimetry,'' 89–102, Contemp. Math. {\bf 545}, Amer. Math. Soc., Providence, RI, 2011

\bibitem[FILS78]{FILS78} Fr\"ohlich, J., R. Israel, E. H. Lieb, and B. Simon, 
{\it Phase transitions and reflection positivity. I. General theory and long range lattice models}, 
Comm. Math. Phys. {\bf 62:1} (1978), 1--34 

\bibitem[FOS83]{FOS83} Fr\"ohlich, J., Osterwalder, K., and E. Seiler,
{\it On virtual representations of symmetric spaces and their analytic
continuation}, Annals Math. {\bf 118} (1983), 461--489


\bibitem[GJ77]{GJ77} Glimm, J., and A. Jaffe, {\it New developments
in quantum field theory and statistical mechanics}, im
``Carg\`ese (1976),'' Plenum, New York, 1977, 35--66

\bibitem[GJ81]{GJ81} --- , ``Quantum Physics--A Functional Integral Point of View,''
Springer-Verlag, New York, 1981

\bibitem[GrN09]{GrN09} Grundling, H., and K.-H.\ Neeb,
{\it Full regularity for a
$C^*$-algebra of the canonical commutation relations},
Reviews in Math. Physics {\bf 21:5} (2009), 587--613


\bibitem[He78]{He78} Helgason, S., ``Differential Geometry, Lie Groups and Symmetric Spaces,'' Academic
Press, 1978.


\bibitem[HN93]{HN93} Hilgert, J., and K.-H. Neeb,
``Lie Semigroups and Their Applications'', Lecture Notes in Mathematics, {\bf 1552}, Springer Verlag, Berlin, 1993.

\bibitem[HN95]{HN95} ---,
{\it Compression semigroups of open orbits
on real flag manifolds}, Monats\-hefte f\"ur Math. {\bf 119} (1995), 187--214

\bibitem[HO96]{HO96} Hilgert, J., and G. \'Olafsson, ``Causal Symmetric Spaces,
Geometry and Harmonic Analysis,'' Acad. Press, 1996

\bibitem[JA08]{JA08} Jaffe, A., {\it Quantum theory and relativity},  in ``Group Representations,
Ergodic Theory, and Mathematical Physics: A Tribute to
George~W. Mackey,'' R. S. Doran, C. C. Moore, R. J. Zimmer, eds.,
Contemp. Math. {\bf 449}, Amer. Math. Soc., 2008.

\bibitem[JR07a]{JR07a} Jaffe, A., and G. Ritter, {\it Quantum field theory 
on curved backgrounds. I. The euclidean functional integral},
Comm. Math. Phys. {\bf 270} (2007), 545--572

\bibitem[JR07]{JR07} ---, {\it Quantum field theory
on curved backgrounds. II. Spacetime symmetries},
arXiv:hep-th/0704.0052v1

\bibitem[JR08]{JR08} ---,  {\it Reflection
positivity and monotonicity},  J.~Math. Phys.  {\bf 49:5}  (2008),
052301, 10 pp.

\bibitem[Jo86]{Jo86} Jorgensen, P. E. T., {\it Analytic continuation of local
representations of Lie groups}, Pac. J. Math. {\bf 125:2} (1986),
397--408

\bibitem[Jo87]{Jo87} ---,  {\it Analytic continuation of local
representations of symmetric spaces}, J.~Funct. Anal. {\bf 70} (1987),
304--322

\bibitem[JO98]{JO98} Jorgensen, P. E. T., and G. \'Olafsson, {\it Unitary representations of Lie groups with
reflection symmetry}, J. Funct. Anal. {\bf 158} (1998), 26--88 

\bibitem[JO00]{JO00} --- , {\it Unitary representations and Osterwalder-Schrader
duality}, in ``The Mathematical Legacy of Harish--Chandra,'' R. S. Doran and V. S. Varadarajan, eds.,
Proc. Symp. in Pure Math. {\bf 68}, Amer. Math. Soc., 2000

\bibitem[Kl77]{Kl77} Klein, A., {\it Gaussian OS-positive processes},
Z. Wahrscheinlichkeitstheorie und Verw. Gebiete {\bf 40:2} (1977), 115--124

\bibitem[Kl78]{Kl78} --- , {\it The semigroup characterization
of Osterwalder--Schrader path spaces and the construction of
euclidean fields}, J. Funct. Anal. {\bf 27} (1978), 277--291

\bibitem[KL81]{KL81} Klein, A., and L. Landau,
{\it Construction of a unique selfadjoint operator for a symmetric
local semigroup}, J. Funct. Anal. {\bf 44} (1981), 121--136

\bibitem[KL82]{KL82} ---, {\it From the Euclidean group to the
Poincar\'e group via Osterwalder-Schrader positivity}, Comm. Math. Phys. {\bf 87} (1982/83),
469--484

\bibitem[KS71]{KS71} Knapp, A. W., and E. M. Stein, 
{\it Intertwining operators for semisimple groups}, Ann. of Math. {\bf 93} (1971), 489--578

\bibitem[Kr49]{Kr49} Krein, M. G., {\it Hermitian positive definite
kernels on homogeneous spaces I, II}, Ukr. Math. J.
{\bf 1} (1949), 64--98;  {\bf 2} (1950), 10--59;
engl. transl. in Amer. Math. Soc. Transl. Ser. 2, Vol. {\bf 34},
69--108, 109--164


\bibitem[Lo69]{Lo69} Loos, O., ``Symmetric Spaces I: General Theory,''
Benjamin, New York, Amsterdam, 1969

\bibitem[LM75]{LM75} L\"uscher, M., and G. Mack, {\it Global conformal 
invariance and quantum field theory}, Comm.  Math. Phys. {\bf 41} (1975), 
203--234 

\bibitem[Mag92]{Mag92} Magyar, M., ``Continuous Linear Representations,''
North-Holland, Mathematical Studies 168, 1992

\bibitem[MN11]{MN11} Merigon, S., and K.-H. Neeb, {\it
Analytic extension techniques for unitary representations of
Banach--Lie groups}, Internat. Math. Research Notices 2012, no. 18, 4260–4300

\bibitem[Ne98]{Ne98} Neeb, K.-H., {\it Operator valued
positive definite kernels on tubes},
Monatshefte f\"ur Math. {\bf 126} (1998), 125--160

\bibitem[Ne00]{Ne00} --- , ``Holomorphy and Convexity in Lie Theory,''
Expositions in Mathematics {\bf 28}, de Gruyter Verlag, Berlin, 2000

\bibitem[Ne10]{Ne10} ---,  {\it
On differentiable vectors for representations of infinite dimensional
Lie groups}, J. Funct. Anal. {\bf 259} (2010), 2814--2855 

\bibitem[Ne11a]{Ne11a} ---, {\it On analytic vectors for unitary
representations of infinite dimensional Lie groups},
Annales de l'Inst. Fourier {\bf 61:5} (2011), 1441--1476 

\bibitem[Ne11b]{Ne11b} ---, {\it Semibounded representations 
of hermitian Lie groups}, Travaux math\'ematiques. {\bf 21} (2012), 29--109

\bibitem[NO98]{NO98} Neeb, K.--H., 
and B. \O{}rsted, {\it Unitary Highest Weight
Representations in Hilbert Spaces of 
Holomorphic Functions on Infinite Dimensional Domains}, J.\ 
Funct.\ Anal. {\bf 156} (1998), 263--300 

\bibitem[OS73]{OS73} Osterwalder, K., and R. Schrader, {\it
Axioms for Euclidean Green's functions.~1},
Comm. Math. Phys. {\bf 31} (1973), 83--112

\bibitem[OS75]{OS75} ---, {\it
Axioms for Euclidean Green's functions. 2},
Comm. Math. Phys. {\bf 31} (1973), 83--112

\bibitem[Sch86]{Sch86} Schrader, R., {\it Reflection positivity for
the complementary series of $\mathrm{SL}(2n,\C)$}, 
Publ. Res. Inst. Math. Sci. {\bf 22} (1986), 119-141.

\bibitem[Scho97]{Scho97} Schottenloher, M., ``A Mathematical Introduction
to Conformal Field Theory,'' Lecture Notes in Physics {\bf m 43},
Springer, 1997

\bibitem[Schw73]{Schw73} Schwartz, L., ``Th\'eorie des distributions,'' 2nd edition, Hermann, Paris, 1973

\bibitem[Sh84]{Sh84} Shucker, D. S., {\it Extensions and generalizations
of a theorem of Widder and the theory of symmetric local semigroups}, J.
Funct. Anal. {\bf 58} (1984), 291--309

\bibitem[SzN70]{SzN70} Sz.-Nagy, B., and C. Foias, ``Harmonic
Analysis of Operators on Hilbert space,'' North-Holland, Amsterdam, London,
1970

\bibitem[Tr67]{Tr67} Treves, F., ``Topological vector spaces, distributions, and
kernels'', Academic Press, New York, 1967

\bibitem[Wa72]{Wa72} Warner, G., ``Harmonic analysis on semisimple Lie groups
I,'' Springer Verlag, Berlin, Heidelberg, New York, 1972


\end{thebibliography}
\end{document}